\documentclass[a4paper,12pt,twoside,notitlepage,leqno]{article}
\usepackage{a4,amsmath,amssymb,setspace,multirow}
\usepackage{graphicx}
\usepackage[linesnumbered,boxed]{algorithm2e}
\def\q0{\theta}
\def\q{\vartheta}

\def\e0{\epsilon}

\def\f0{\phi}
\def\f{\varphi}

\def\R{{\mathbb R}}
\def\Z{{\mathbb Z}}

\def\e{\varepsilon}

\newtheorem{thm}{Theorem}[section]
\newtheorem{cor}[thm]{Corollary}
\newtheorem{lem}[thm]{Lemma}
\newtheorem{prop}[thm]{Proposition}
\newtheorem{defn}[thm]{Definition}
\newtheorem{rem}[thm]{Remark}

\newcommand\blfootnote[1]{%
  \begingroup
  \renewcommand\thefootnote{}\footnote{#1}%
  \addtocounter{footnote}{-1}%
  \endgroup
}
\begin{document}
\title{\sc The moving frame on the fractal curves\footnote{This work was supported by NSFC(Nos 11201056 and 11371080).}}
\author{Yun Yang, Yanhua Yu\footnote{Corresponding author.}\blfootnote{E-mail addresses: yangyun@mail.neu.edu.cn (Y. Yang), yuyanhua@mail.neu.edu.cn(Y. Yu).}
\\{\small Department of Mathematics, Northeastern University, Shenyang 110004, P. R. China}}
\markboth{}{}
\date{}
\maketitle
\numberwithin{equation}{section}
\begin{abstract}
Using the moving frame and invariants, any discrete curve in $\R^3$ could be uniquely identified by its centroaffine curvatures and torsions\cite{Y-Y}. In this paper, depending on the affine curvatures of the fractal curves, such as Koch curve and Hilbert curve, we can clearly describe their iterative regularities. Interestingly, by the affine curvatures, the fractal curves can be quantified and encoded accordingly to a sequence. Hence, it is more convenient for future reference. Given three starting points, we can directly generate the affine Koch curve and affine Hilbert curve at the step $n, \forall n\in \mathbb{Z}^+$. Certainly,  if the initial three points are standard, the curve is the traditional Koch curve or Hilbert curve. By this method, the characteristic of some fractal curves which look like irregular could be quantified, and the regularities would become more obvious.
\end{abstract}

\section{Introduction}
Discrete differential geometry has attracted much attention recently, mainly due to the growth of computer graphics. One of the main issues in discrete differential geometry is to define suitable discrete analogous of the concepts of smooth differential geometry\cite{Bobenko-2,Bobenko-4}. More recently, the expansion of computer graphics and applications in mathematical physics have given a great impulse to the issue of giving discrete equivalents of affine differential geometric objects\cite{Bobenko-1,Craizer-1,Craizer-2}.
\par Group based moving frames have a wide range of applications, from the classical equivalence problems in differential geometry to more modern applications such as computer vision\cite{Mansfield,Olver}. The first results for the computation of discrete invariants using group based moving frames were given by Olver\cite{Olver-0} who calls them joint invariants; modern applications to date include computer vision\cite{Olver-01} and numerical schemes for systems with a Lie symmetry\cite{Chhay}. Moving frames for discrete applications as formulated by Olver do give generating sets of discrete invariants, and the recursion formulas for differential invariants are so successful for the application of moving frames to calculus-based applications. Recent development of a theory of discrete equivariant moving frames has been
applied to integrable differential-difference systems\cite{Mansfield,Olver}.
\par Following the ideas of Klein, presented in his famous lecture at Erlangen, several geometers in the early 20th century proposed the study
of curves and surfaces with respect to different transformation groups. In geometry, an affine transformation, affine map or an affinity is a function between affine spaces which preserves points, straight lines and planes. Also, sets of parallel lines remain parallel after an affine transformation. An affine transformation does not necessarily preserve angles between lines or distances between points, though it does preserve ratios of distances between points lying on a straight line. Examples of affine transformations include translation, scaling, homothety, similarity transformation, reflection, rotation, shear mapping, and compositions of them in any combination and sequence.
\par A fractal is a never-ending pattern. Fractals are infinitely complex patterns that are self-similar across different scales. They are created by repeating a simple process over and over in an ongoing feedback loop. Driven by recursion, fractals are images of dynamic systems $-$ the pictures of Chaos. Geometrically, they exist in between our familiar dimensions. Fractal patterns are extremely familiar, since nature is full of fractals. For instance: trees, rivers, coastlines, mountains, clouds, seashells, hurricanes, etc. Abstract fractals, such as the Mandelbrot Set, can be generated by a computer calculating a simple equation over and over.
\par The category of fractal curves concerns curves that are fractal in some way, typically by being non-differentiable everywhere, or by having a period-doubling symmetry.
Most of fractal curves are discrete polygonal lines, such as Koch curve, Hilbert curve, dragon curve, etc. Exactly, we have to admit that it is difficult to study those discrete fractal curves which look like irregular. They probably are equivalent to the known fractal curves. Now, using the moving frame and invariants, any discrete curve in $\R^3$ could be uniquely identified by its centroaffine curvatures and torsions\cite{Y-Y}. Thus, the discrete fractal curves can be quantified, and the regularities are more obvious.
\par The arrangement of the paper is as follows: In Sect. 2 we recall the basic theory and notions for affine geometry, the basic notations for discrete curves and centroaffine curves. There are some results with centroaffine curvatures and torsions for centroaffine planar curves and space curves. In Sect. 3, we consider the affine Koch curve and its regularities. There are some examples for the affine Koch curves. In Sects. 4, we study the affine Koch snowflake and obtain its iterative sequences.  Finally, in Sect. 5 we describes the affine Hilbert curve, and then obtain the sequence for every affine Hilbert curve.
\section{Preliminaries}
\subsection{Affine mappings and transformation groups, basic notations}
If $X$ and $Y$ are affine spaces, then every affine transformation $f:X\rightarrow Y$  is of the form $\vec{x}\mapsto M\vec{x}+\vec{b}$ , where $M$ is a linear transformation on $X$ and  $b$ is a vector in $Y$. Unlike a purely linear transformation, an affine map need not preserve the zero point in a linear space. Thus, every linear transformation is affine, but not every affine transformation is linear.
\par For many purposes an affine space can be thought of as Euclidean space, though the concept of affine space is far more general (i.e., all Euclidean spaces are affine, but there are affine spaces that are non-Euclidean). In affine coordinates, which include Cartesian coordinates in Euclidean spaces, each output coordinate of an affine map is a linear function (in the sense of calculus) of all input coordinates. Another way to deal with affine transformations systematically is to select a point as the origin; then, any affine transformation is equivalent to a linear transformation (of position vectors) followed by a translation.
\par It is well known that the set of all automorphisms of a vector space $V$ of dimension $m$ forms a group. We use the following standard notations for this group and its subgroups(\cite{L-U-Z}):
$$GL(m,\R):=\{L:V\rightarrow V|L\quad isomorphism\};$$
$$SL(m,\R):=\{L\in GL(m,\R)|\det L=1\}.$$
Correspondingly, for an affine space $A, \dim A=m$, we have the following affine transformation groups.
\begin{flalign*}
  &\mathcal{A}(m):=\{\alpha:A\rightarrow A|L_{\alpha}\ \mathrm{regular}\}\quad \mathrm{is\ the\ regular\ affine\ group}.  \\
  &\mathcal{S}(m):=\{\alpha\in \mathcal{A}|\det\alpha=1\}\quad \mathrm{is\ the\ unimodular(equiaffine)\ group}.  \\
  &\mathcal{Z}_p(m):=\{\alpha\in \mathcal{A}|\alpha(p)=p\}\quad \mathrm{is\ the\ centroaffine\ group\ with\ center}\ p\in A. \\
  &\tau(m):=\{\alpha:A\rightarrow A| \mathrm{there\ exists}\ b(\alpha)\in V, \mathrm{s.t.}\ \overrightarrow{p\alpha(p)}=b(\alpha), \forall p\in A\}\\
  &\qquad\qquad \mathrm{is\ the\ group\ of\ transformations\ on\ } A.
\end{flalign*}
Let $\mathcal{G}$ be one of the groups above and $S_1,S_2\subset A$ subsets. Then $S_1$ and $S_2$ are called equivalent modulo $\mathcal{G}$ if there exists an $\alpha\in\mathcal{G}$ such that $$S_2=\alpha S_1.$$
\par In centroaffine geometry we fix a point in $A$(the origin $O\in A$ without loss of generality) and consider the geometric properties in variant under the centroaffine group $\mathcal{Z}_p$. Thus the mapping $\pi_0:A\rightarrow V$ identifies $A$ with the vector space $V$ and $\mathcal{Z}_O$ with $GL(m, \R)$.
\subsection{Discrete curves and the affine curvatures}
The aim of this section is to collect some useful definitions and results(for a more detailed exposition of this material, please refer to \cite{Y-Y}).
We will start very simple, by discretizing the notion of a smooth curve. That is, we want to define a discrete analog to a smooth map from an interval
$I\subset \R$ to $\R^n$. By discrete we mean here that the map should not be defined on an interval in $\R$ but on a discrete (ordered) set of points therein. It
turns out that this is basically all we need to demand in this case:
\begin{defn} Let $I\subset\Z$ be an interval (the intersection of an interval in $\R$ with $\Z$, possibly infinite). A map $\vec{r}:I\rightarrow \R^n$ is called a discrete curve, when we put the starting point of the vector $\vec{r}$ to the origin $O\in\R^n$. Obviously, a discrete curve is  a polygon. A discrete curve $\vec{r}$ is said to be periodic (or closed) if $I=\Z$ and if there is a $p\in\Z$ such that $\vec{r}(k)=\vec{r}(k+p)$ for all $k\in I$. The smallest possible value of $p$ is called the period.
\end{defn}
In fact, we can define
\begin{equation}
  \vec{r}(t)=(t-k)\vec{r}(k)+(k+1-t)\vec{r}(k+1), \quad\forall t\in (k,k+1), k\in\Z.
\end{equation}
Then simplicity of a smooth curve can be generalized to the discrete case.
We shall need a definition of the discrete centroaffine curves in $\R^2$ and $\R^3$.
\begin{defn}\label{def-CA}
A discrete planar curve $\vec{r}:I\rightarrow \R^2$ is called a centroaffine planar curve if the edge tangent vector $\vec{t}_k$  is not parallel to position vectors $\vec{r}(k)$ and $\vec{r}(k+1)$, and a discrete curve $\vec{r}:I\rightarrow \R^3$ is called a centroaffine curve if the edge tangent vectors $\vec{t}_{k-1}, \vec{t}_k$ and the position vector $\vec{r}(k)$ are not coplanar.
\end{defn}
In \cite{Y-Y}, we consider two invariants $\kappa, \bar{\kappa}$ and their geometrical properties under the affine transformation, although we call them the first and second centroaffine curvatures just for unity, that because for the space discrete curve, we will use them together with centroaffine torsions, only under the centroaffine transformation. In the next section, we directly call them affine curvatures if the discrete curve is planar. Now let vector-valued function $\vec{r}:I\subset\Z\rightarrow \R^2$ represent a discrete planar curve $C$.
\begin{defn}For the discrete planar curve $C$,  if $[\vec{t}_{k-1},\vec{t}_{k}]=0$, we call its first centroaffine curvature $\kappa_k=0$ at point $\vec{r}(k)$, which implies the curve is a straight line locally to $\vec{r}(k)$, where $[\cdots]$ denotes the standard determinant in $\R^2$. If $[\vec{t}_{k-1},\vec{t}_{k}]\neq0$,  the first and second centroaffine curvatures at the point $\vec{r}(k)$ are defined by
\begin{equation}\label{PCur}
  \kappa_k=\frac{[\vec{t}_{k},\vec{t}_{k+1}]}{[\vec{t}_{k-1},\vec{t}_{k}]},\quad \bar{\kappa}_k=\frac{[\vec{t}_{k-1},\vec{t}_{k+1}]}{[\vec{t}_{k-1},\vec{t}_{k}]}.
\end{equation}
\end{defn}
 \par In fact, from Eq. (\ref{PCur}), we can obtain the chain structure
 \begin{equation}\label{Ite}
   \vec{r}_{k+2}-\vec{r}_{k+1}=-\kappa_k(\vec{r}_k-\vec{r}_{k-1})+\bar{\kappa}_k(\vec{r}_{k+1}-\vec{r}_k).
 \end{equation}
 This shows that
 \begin{equation}\label{Ite-Pt}
   \vec{r}_{k+2}=\kappa_k\vec{r}_{k-1}+(-\kappa_k-\bar{\kappa}_k)\vec{r}_k+(1+\bar{\kappa}_k)\vec{r}_{k+1}.
 \end{equation}
Let curve $\vec{r}:I\subset\Z\rightarrow \R^3$ be a centroaffine discrete curve denoted by $C$, and then by the definition \ref{def-CA}, we have $[\vec{r}_k, \vec{t}_{k-1},\vec{t}_{k}]\neq0$, where $[\cdots]$ denotes the standard determinant in $\R^3$. In the following
 the centroaffine curvatures and centroaffine torsions of a centroaffine discrete space curve in $\R^3$ will be defined.
\begin{defn} The first, second centroaffine curvatures and centroaffine torsions of the discrete cnetroaffine curve $\vec{r}$ at point $\vec{r}(k)$ are defined by
\begin{equation}\label{PCT}
  \kappa_k:=\frac{[\vec{r}_{k+1},\vec{t}_{k},\vec{t}_{k+1}]}{[\vec{r}_k,\vec{t}_{k-1},\vec{t}_{k}]},\quad \bar{\kappa}_k:=\frac{[\vec{r}_{k+1},\vec{t}_{k-1},\vec{t}_{k+1}]}{[\vec{r}_k, \vec{t}_{k-1},\vec{t}_{k}]},\quad \tau_k:=\frac{[\vec{t}_{k-1},\vec{t}_{k},\vec{t}_{k+1}]}{[\vec{r}_k, \vec{t}_{k-1},\vec{t}_{k}]}.
\end{equation}
\end{defn}
By Eq. \ref{PCT}, under a centroaffine transformation $\R^3\ni \vec{x}\mapsto A\vec{x}\in\R^3$, where $A\in GL(3,\R)$, it is easy to see that the first, second centroaffine curvatures and centroaffine torsions are invariant.
However, under an affine transformation $\vec{x}\mapsto A\vec{x}+\vec{b}$, where $\vec{b}\in\R^3$ is a constant vector, the first, second centroaffine curvatures and centroaffine torsions may change. Hence, we have
\begin{prop}
The first, second centroaffine curvatures and centroaffine torsions are centroaffine invariants and not affine invariants\cite{Y-Y}.
\end{prop}
In fact,
$$ [\vec{r}_k,\vec{t}_{k-1},\vec{t}_{k}]=[\vec{r}_{k-1}, \vec{r}_{k},\vec{r}_{k+1}].$$
Then the centroaffine curvatures and torsions can be rewritten as
\begin{equation}\label{PCT-1}
  \kappa_k=\frac{[\vec{r}_{k},\vec{r}_{k+1},\vec{r}_{k+2}]}{[\vec{r}_{k-1},\vec{r}_{k},\vec{r}_{k+1}]},\quad \bar{\kappa}_k=\frac{[\vec{r}_{k+1},\vec{t}_{k-1},\vec{r}_{k+2}]}{[\vec{r}_{k-1},\vec{r}_{k},\vec{r}_{k+1}]},\quad
  \tau_k=\frac{[\vec{t}_{k-1},\vec{t}_{k},\vec{t}_{k+1}]}{[\vec{r}_{k-1}, \vec{r}_{k},\vec{r}_{k+1}]}.
\end{equation}
By a direct calculation, it follows that
\begin{equation}\label{Ite-3d-P}
  \vec{r}_{k+2}=\kappa_k\vec{r}_{k-1}+(-\kappa_k-\bar{\kappa}_k)\vec{r}_{k}+(\tau_k+\bar{\kappa}_k+1)\vec{r}_{k+1}, \quad \forall k\in \Z,
\end{equation}
and
\begin{equation}\label{Ite-3d-MP}
  (\vec{r}_{k+2}, \vec{r}_{k+1}, \vec{r}_{k})=  (\vec{r}_{k+1}, \vec{r}_{k}, \vec{r}_{k-1})\left(
                                                                                             \begin{array}{ccc}
                                                                                               \tau_k+1+\bar{\kappa}_k & 1 & 0 \\
                                                                                               -\kappa_k-\bar{\kappa}_k & 0 & 1 \\
                                                                                               \kappa_k & 0 & 0 \\
                                                                                             \end{array}
                                                                                           \right),
\end{equation}
which are the three dimensional curve chain structures. This formula is called the Frenet-Serret formula of a discrete centroaffine curve.
\par On the other hand, when $\tau=0$, from Eq. (\ref{Ite-3d-P}) we get the chain of edge tangent vector
\begin{equation}\label{Ite-3d-PT}
  \vec{t}_{k+1}=-\kappa_k\vec{t}_{k-1}+\bar{\kappa}_k\vec{t}_{k}, \quad \forall k\in \Z,
\end{equation}
which is coincident with Eq. (\ref{Ite}).
\par If $\kappa_k\neq0$, we notice that the inverse chain can be represented as
\begin{equation}\label{Ite-3d-MPI}
  (\vec{r}_{k+1}, \vec{r}_{k}, \vec{r}_{k-1})=  (\vec{r}_{k+2}, \vec{r}_{k+1}, \vec{r}_{k})\left(
                                                                                             \begin{array}{ccc}
                                                                                               0 & 0 & \frac{1}{\kappa_k} \\
                                                                                               1 & 0 & -\frac{\tau_k+1+\bar{\kappa}_k}{\kappa_k} \\
                                                                                               0 & 1 & 1+\frac{\bar{\kappa}_k}{\kappa_k} \\
                                                                                             \end{array}
                                                                                           \right).
\end{equation}
In this section, we only consider the discrete centroaffine space curve $C$ under the centroaffine transformation. Firstly, the following proposition states that with the given first, second centroaffine curvatures and centroaffine torsions, the curve is only determined up to a centroaffine transformation.
\begin{prop}\label{space-ceneq}
Two curve $C$ and $\bar{C}$ are centroaffine equivalent if and only if they have same centroaffine curvatures $\kappa_k, \bar{\kappa_k}$ and torsions $\tau_k$,
for all $k\in I\subset\Z$\cite{Y-Y}.
\end{prop}
\begin{rem}\label{rem-planar}
If $\tau_k=0, \forall k\in\Z$, the curve $C$ is a planar curve and the centroaffine curvatures in Eq. (\ref{PCT}) are as same as defined in Eq. (\ref{PCur}) for the planar curve. Hence, if the discrete curve is planar, we directly call them affine curvatures, since they are invariant under the affine transformation.
\end{rem}
\section{Koch curves}
\par 	The Von Koch curves, named from the Swedish mathematician Helge Von Koch who originally devised them in 1904, are perhaps the most beautiful fractal curves. These curves are amongst the most important objects used by Benoit Mandelbrot for his pioneering work on fractals. More than any other, the Von Koch curves allows numerous variations and have inspired many artists that produced amazing pieces of art.
\par
\begin{figure}[hbtp]
            \centering
            \begin{tabular}{ccc}
              \includegraphics[width=.3\textwidth]{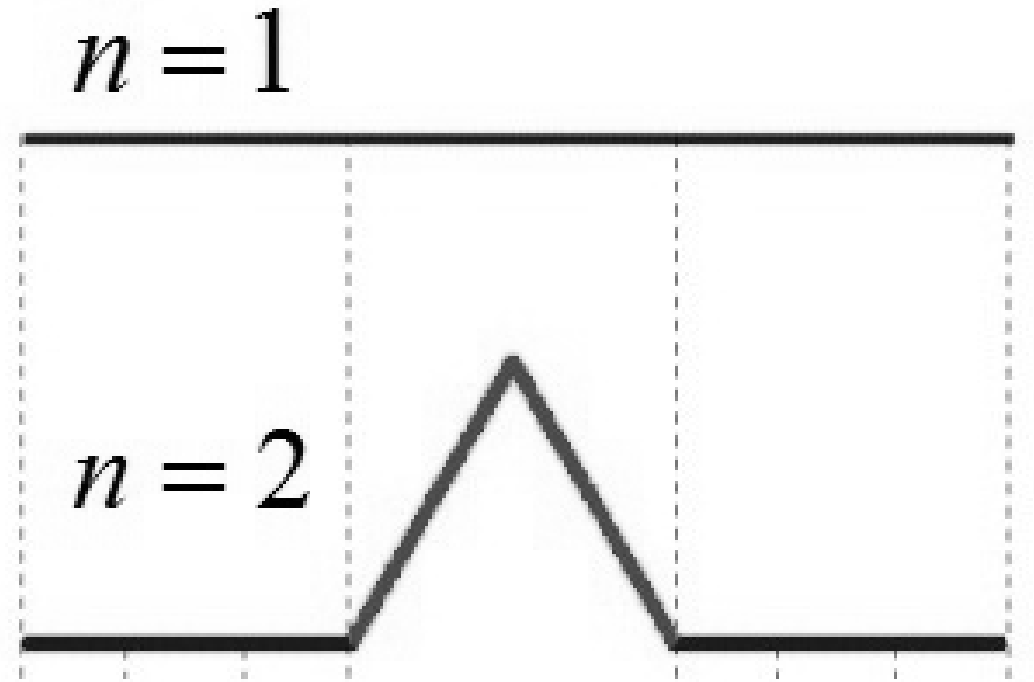}&\includegraphics[width=.3\textwidth]{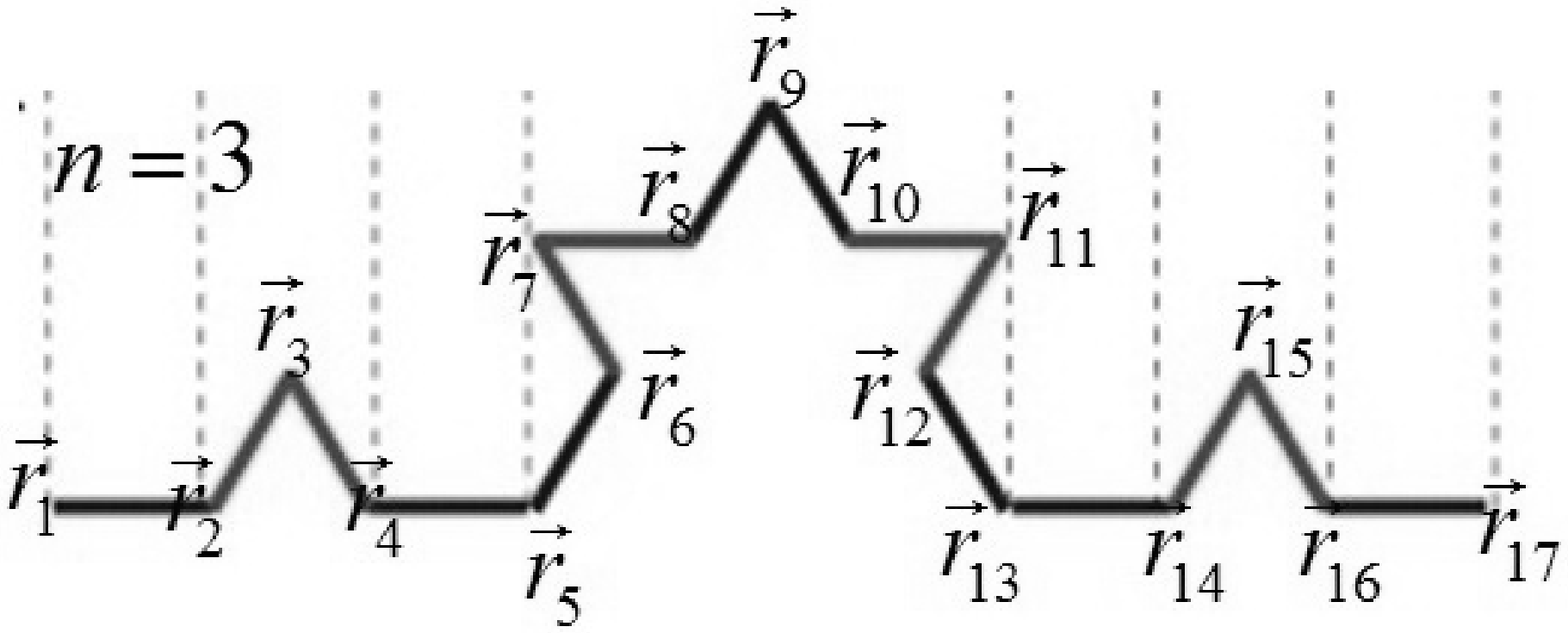}& \includegraphics[width=.25\textwidth]{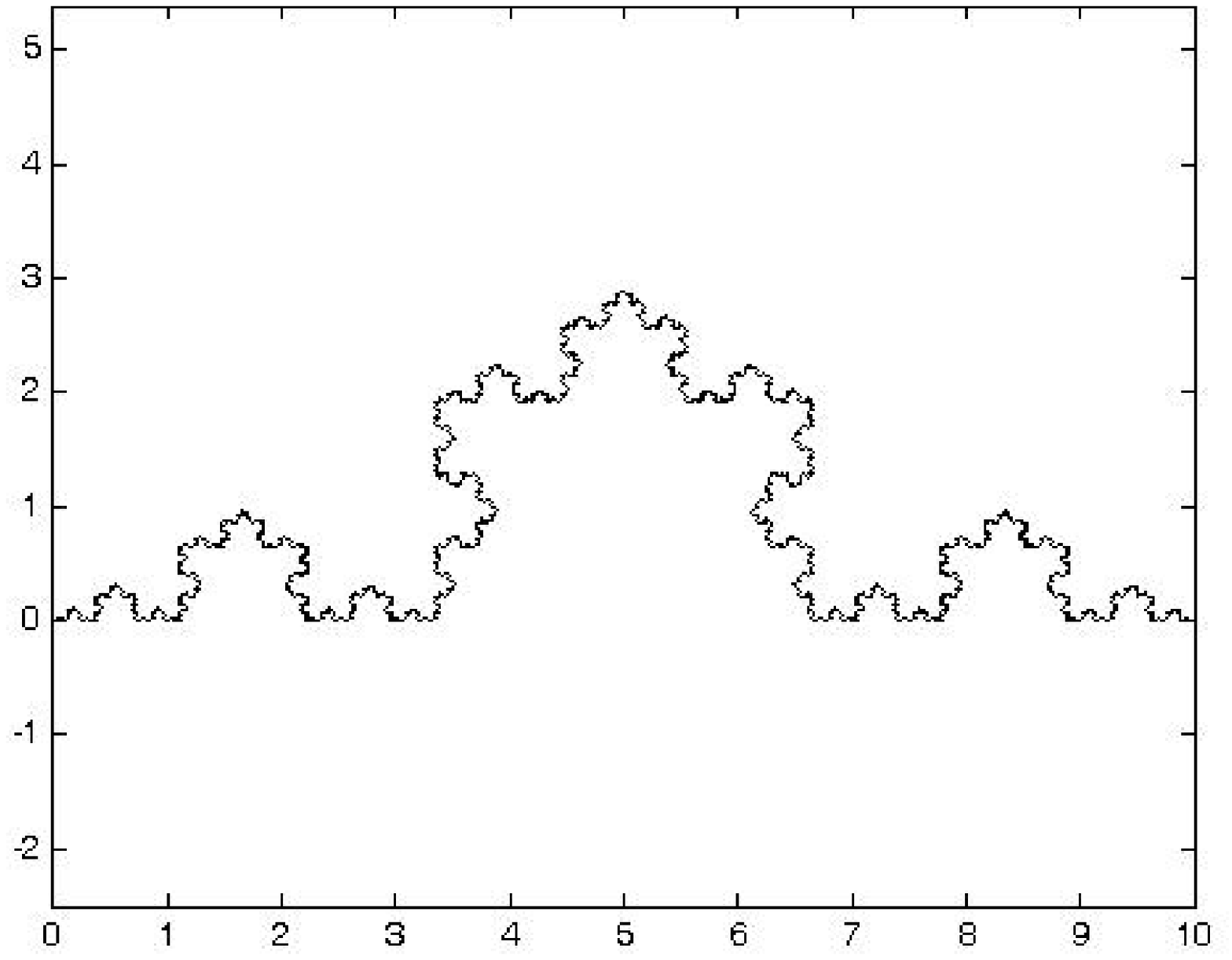}
            \end{tabular}
            \caption{Construction of the Koch curve. }
            \label{fig-Koch}
 \end{figure}

 The construction of the curve is fairly simple. A straight line is first divided into three equal segments. The middle segment is removed and replaced by two segments having the same length to generate an equilateral triangle. Now repeat, taking each of the four resulting segments, dividing them into three equal parts and replacing each of the middle segments by two sides of an equilateral triangle. Continue this construction. The process is shown in Figure \ref{fig-Koch}.
\par
\begin{figure}[hbtp]
            \centering
            \begin{tabular}{cc}
              \includegraphics[width=.35\textwidth]{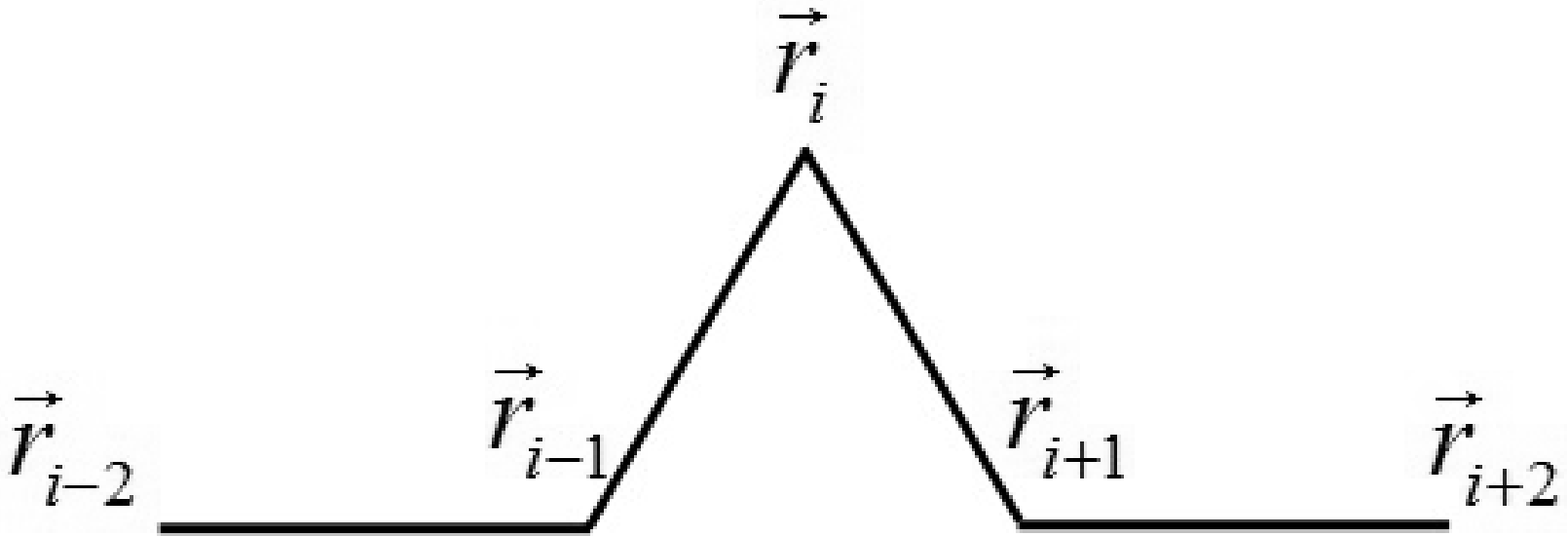}\quad& \includegraphics[width=.2\textwidth]{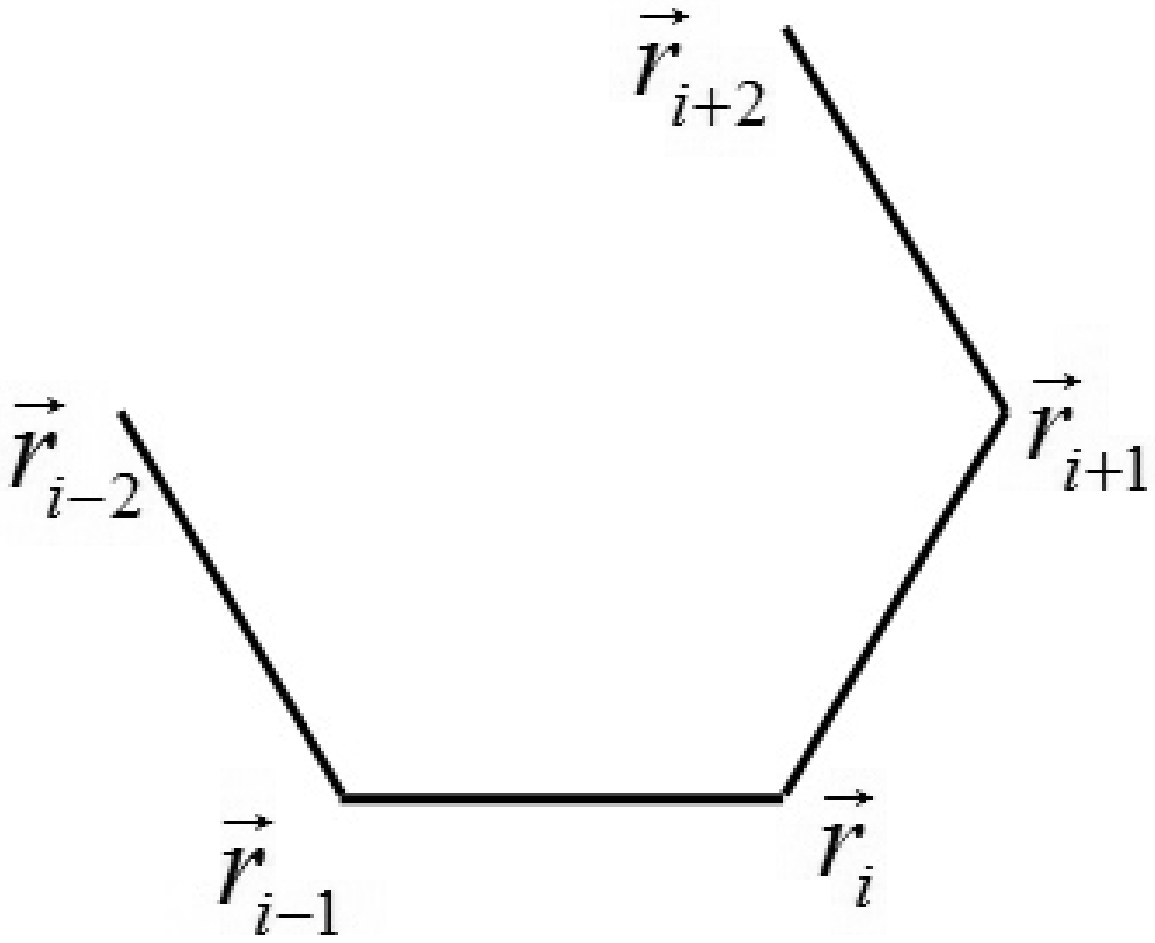}
            \end{tabular}
            \caption{{\it Left:} Sharp corner in Koch curve and sharp angle pair $\{\vec{r}_{i-1}, \vec{r}_i\}$. {\it Right:} Obtuse corner in Koch curve and obtuse angle pair $\{\vec{r}_{i-1}, \vec{r}_i\}$. }
            \label{fig-corner}
 \end{figure}

 Obviously, Koch curve is comprised of a large number of sharp corners and obtuse corners as shown in Figure \ref{fig-corner}. Now we call the pair $\{\vec{r}_{i-1}, \vec{r}_i\}$ shown in the left of Figure \ref{fig-corner} {\it sharp angle pair}, and the pair $\{\vec{r}_{i-1}, \vec{r}_i\}$ shown in the right of Figure \ref{fig-corner} {\it obtuse angle pair}. The point $\vec{r}_{i}$ in a sharp angle pair $\{\vec{r}_{i-1}, \vec{r}_i\}$ is called {\it sharp point}.
 \par For example, in the middle of Figure \ref{fig-Koch}, the sharp angle pairs include $\{\vec{r}_{2}, \vec{r}_3\}$, $\{\vec{r}_{6}, \vec{r}_7\}$, $\{\vec{r}_{8}, \vec{r}_9\}$, $\{\vec{r}_{10}, \vec{r}_{11}\}$ and $\{\vec{r}_{14}, \vec{r}_{15}\}$. Obviously, $\{\vec{r}_{4}, \vec{r}_5\}$ and $\{\vec{r}_{12}, \vec{r}_{13}\}$ are obtuse angle pairs. There are three marginal points $\vec{r}_1, \vec{r}_{16}, \vec{r}_{17}$ and five sharp points $\vec{r}_3,\vec{r}_7, \vec{r}_9, \vec{r}_{11}, \vec{r}_{15}$.
 \par Using Eq. (\ref{PCur}), by a direct computation, for a sharp angle pair $\{\vec{r}_{i-1}, \vec{r}_i\}$, we have
 \begin{equation}\label{KCVC}
   \kappa_{i-1}=\kappa_{i}=-1, \quad \bar{\kappa}_{i-1}=-1, \quad \bar{\kappa}_i=1.
 \end{equation}
 Similarly, for a obtuse angle pair $\{\vec{r}_{i-1}, \vec{r}_i\}$, we can obtain
 \begin{equation}\label{KCCC}
   \kappa_{i-1}=\kappa_{i}=1, \quad \bar{\kappa}_{i-1}=\bar{\kappa}_i=1.
 \end{equation}
As we know, at the step $n$, there are $4^{n-1}+1$ points, where $\vec{r}_1, \vec{r}_{4^{n-1}}$ and $\vec{r}_{4^{n-1}+1}$ are marginal points. There is not the affine curvature at these marginal points. Now it is not difficult to list the sequence numbers of sharp points at the step $n$, which are
\begin{equation}\label{KCP}
  4^{n-2-i}\times(4j-2)+1, \quad i=0,1,2,\cdots, n-2, j=1,2,3,\cdots,4^i.
\end{equation}
Hence we conclude that
\begin{prop}
For a Koch curve with $4^{n-1}+1$ points at the step $n$, there are $\frac{4^{n-1}-1}{3}$ sharp points. Therefore, it has $\frac{4^{n-1}-1}{3}$ sharp angle pairs and $\frac{4^{n-1}-4}{6}$ obtuse angle pairs, where $n\geq 2$.  The sequence numbers of sharp points are listed in Eq. (\ref{KCP}).
\end{prop}
Obviously,
\begin{equation}
  \lim_{n\rightarrow \infty}\frac{\frac{4^{n-1}-1}{3}}{\frac{4^{n-1}-4}{6}}=2.
\end{equation}
Since the curvatures of the sharp angle pair and the obtuse angle pair are certain,  we can treat the two points $\vec{r}_{i-1}, \vec{r}_i$ of a sharp angle pair as one element, which is marked by $1$. Similarly, we mark the two points $\vec{r}_{i-1}, \vec{r}_i$ of a obtuse angle pair as $0$. The first point $\vec{r}_1$ and the last two points $\vec{r}_{4^{n-1}}, \vec{r}_{4^{n-1}+1}$ are neglected because they do not have the affine curvatures.
\par Then, $n=2$, we get the affine Koch curve denoted as $1$. When $n=3$, the affine Koch curve is denoted by $1011101$. Exactly, $n=4$, the affine Koch curve can be represented  by
$$1011101{\bf0}1011101{\bf1}1011101{\bf0}1011101.$$
Obviously, we could find the rule of the iteration. That is, if we denote the affine Koch curve of the step $n$ by $A$, which includes $2\times4^{n-2}-1$ elements, at the step $n+1$, the affine Koch curve can be written as $A0A1A0A$. So we obtain the following results.
\begin{prop}\label{prop-Ko}
An affine Koch curve can be encoded as one of the following forms $$1, 1{\bf 0}1{\bf1}1{\bf0}1, 1011101{\bf0}1011101{\bf1}1011101{\bf0}1011101, \cdots\cdots,$$
where the element $1$ means there are two points with affine curvatures $\kappa_1=\kappa_2=-1, \bar{\kappa}_1=-1, \bar{\kappa}_2=1$, and the element $0$ implies the two points have the affine curvatures $\kappa_1=\kappa_2=\bar{\kappa}_1=\bar{\kappa}_2=1$. Further more, from the bold number we can find the characterization of the iteration. There are $2\times4^{n-2}-1$ elements at the step $n$, where $n\ge 2$.
\end{prop}
By the above results, we have
\begin{prop}
After one iteration, an affine Koch curve will generate 4 times of itself, plus two new obtuse angel pairs and one new sharp angle pair.
\end{prop}
Now, Let us find the positions of the element $1$ in the sequence of Proposition \ref{prop-Ko}. From Eq. (\ref{KCP}), it is not difficult to see that at the position
\begin{equation}\label{Koch-idx}
idx=4^i(2k-1), i=0,1,2,\cdots, n-2, k=1,2,\cdots ,4^{n-2-i},
\end{equation}
the element is $1$, where $n\ge2$. So it is easy to see
\begin{cor}At the sequences of the affine Koch curve, the position of the element $1$ satisfies Eq. (\ref{Koch-idx}). Hence, there are $\frac{4^{n-1}-1}{3}$ element $1$'s and $\frac{2(4^{n-2}-1)}{3}$ element $0$'s at the $n^{\mathrm{th}}$ sequence, where $n\ge2$.
\end{cor}
Hence, using the affine curvatures, we obtain the algorithm \ref{alg-koch}, which can directly generate an affine Koch curve only depending on the initial three starting points for any arbitrary step $n$.
\begin{algorithm}[hbtp]
\label{alg-koch}
 \KwIn{initial three starting points $P_1, P_2, P_3$, and an integer $N>1$;}
 \KwOut{an affinely Koch curve of the $N$th step;}
 join points $P_1$, $P_2$ and $P_3$ with two straight line segments in turn\;
 \For{$idx\Leftarrow1$ \KwTo $2\times4^{N-1}-1$}{
      $i\Leftarrow 0, flag\Leftarrow 0$\;
     \While{$i<=\log_4 idx$ {\bf and} $flag==0$}{
         \If{$\frac{idx}{4^i}$ is odd}
         {$flag\Leftarrow 1$;
         }
     $i\Leftarrow i+1$;
     }
     \eIf{$flag==1$}
     {use the affine curvatures $\kappa=\{-1,-1\},\bar{\kappa}=\{-1,1\} $ and Eq. (\ref{Ite-Pt}) to generate two points, start from $P_3$, join them with the straight line segments in turn, assign into $P_1,P_2,P_3$ with the last three points\;
     }{use the affine curvatures $\kappa=\{1,1\},\bar{\kappa}=\{1,1\} $ and Eq. (\ref{Ite-Pt}) to generate two points, start from $P_3$, join them with the line segments in turn, assign into $P_1,P_2,P_3$ with the last three points\;
     }
  }
\caption{Generate directly an affine Koch curve by the affine curvatures}
\end{algorithm}

\begin{figure}[hbtp]
            \centering
            \begin{tabular}{cc}
              \includegraphics[width=.26\textwidth]{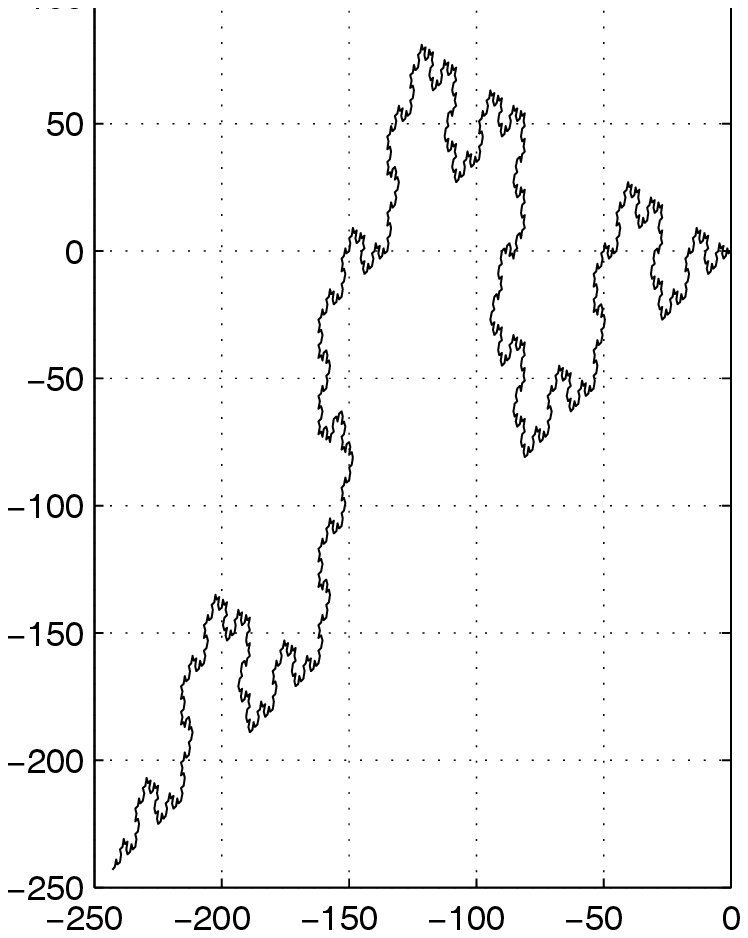}\quad& \includegraphics[width=.55\textwidth]{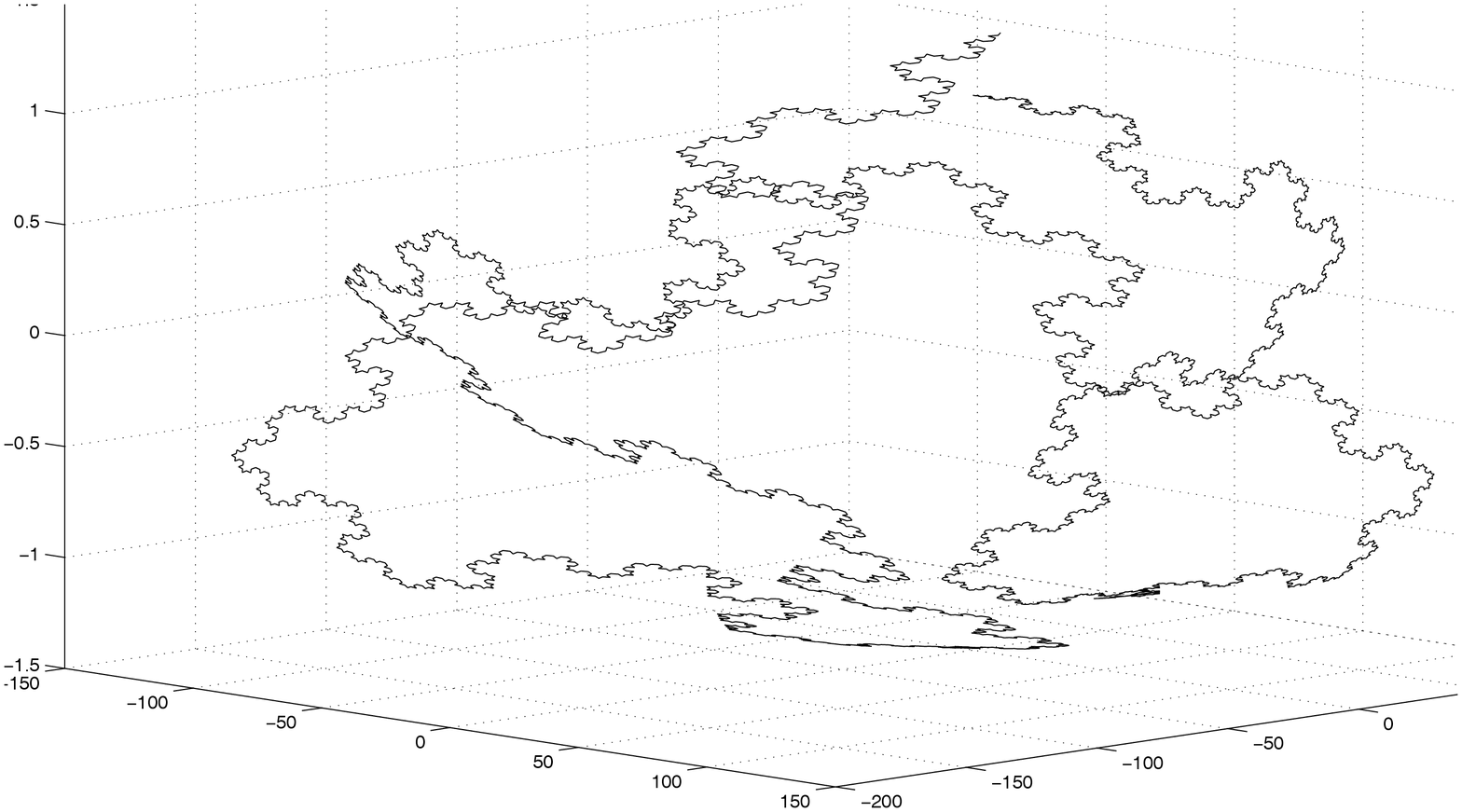}\\ \includegraphics[width=.4\textwidth]{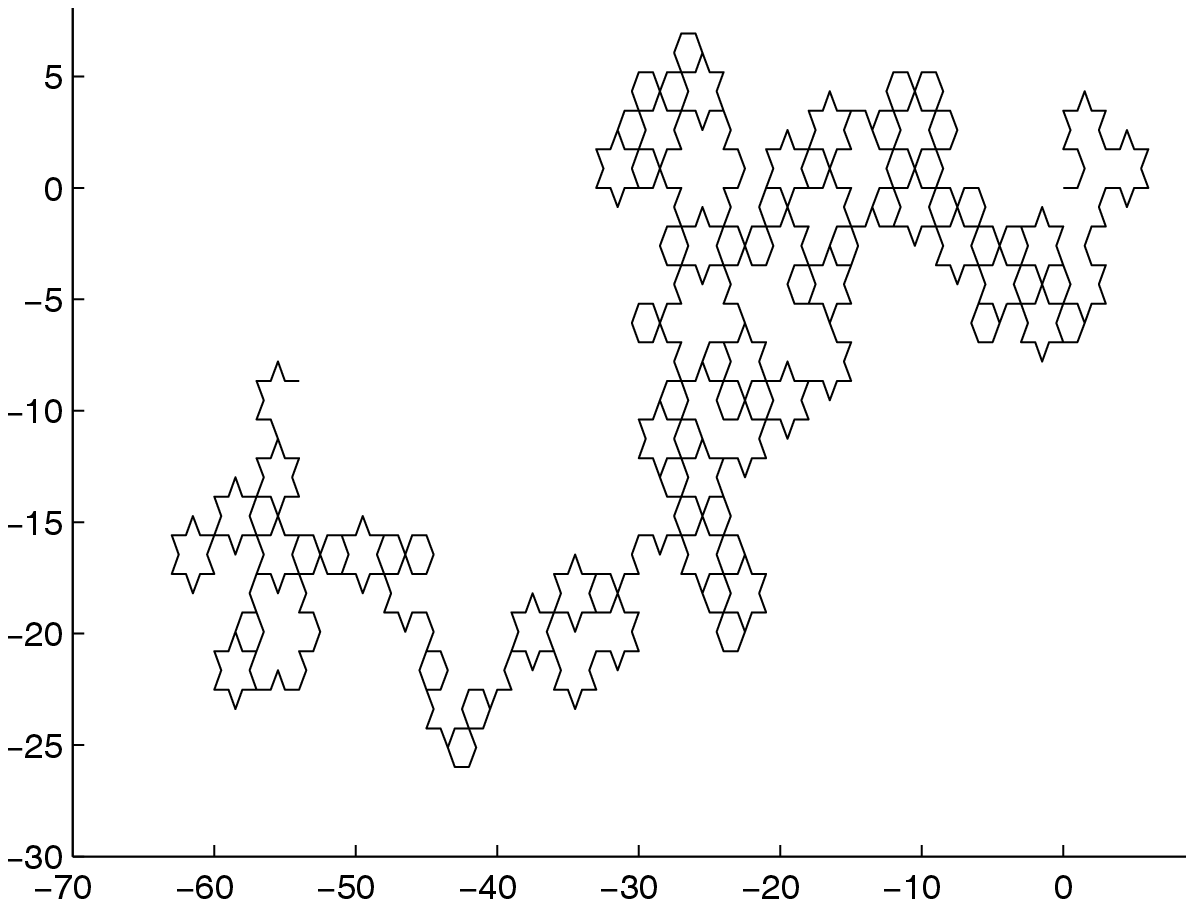}\quad&\includegraphics[width=.45\textwidth]{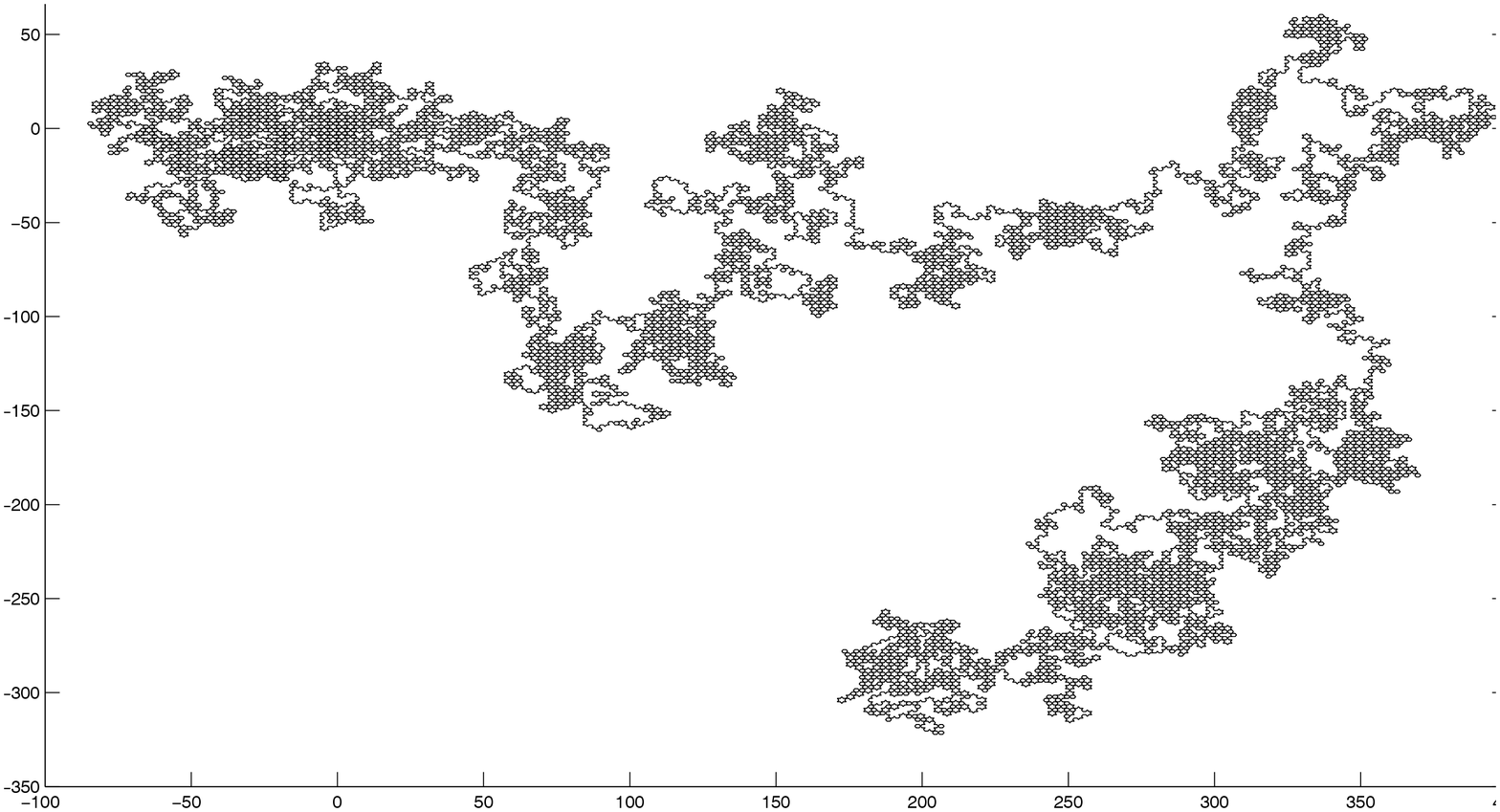}
            \end{tabular}
            \caption{Discrete affine curve generated by the affine curvatures of Koch curve. }
            \label{fig-affineK}
 \end{figure}

\begin{rem}
Now using the algorithm \ref{alg-koch}, we could get an affine Koch curve. If we choose three initial points $\vec{r}_1,\vec{r}_2$ and $\vec{r}_3$ which are not collinear, the curve can be determined uniquely. If we choose $\vec{r}_3=\vec{r}_2+\left(
                       \begin{array}{cc}
                         \cos\frac{\pi}{3} & -\sin\frac{\pi}{3} \\
                         \sin\frac{\pi}{3} & \cos\frac{\pi}{3} \\
                       \end{array}
                     \right)
(\vec{r}_2-\vec{r}_1),$ it is a standard Koch curve. Otherwise, it is affinely equivalent to the standard Koch curve.
\end{rem}
The following example shows how to use the affine curvatures of the Koch curve to generate fractal curves.\\
\noindent {\bf Example 1.} In this example, we give four curves which are obtained by the affine curvatures of the affine Koch curve:
\begin{description}
  \item[(a)] the first graph of Figure \ref{fig-affineK} is generated by using the initial points $[0~0;-1~-1;-1.5~1]$ and $n=5$, which is an affine Koch curve;
  \item[(b)] using the the initial points $[0~0~1;1~0~1;1.5~\frac{\sqrt{3}}{2}~1]$ and $n=7$, with the same affine curvatures as the affine Koch curve, and the centroaffine torsions $\tau_i=0.0001\kappa_i$, we get a discrete space curve as shown in the second graph of Figure \ref{fig-affineK};
  \item[(c)] the last two graphs in Figure \ref{fig-affineK} show two stochastic affine Koch curves generated by the same proportion sharp points as the affine Koch curve with $n=5$ and $n=8$, where we use the initial points $[0~0;1~0;1.5~\frac{\sqrt{3}}{2}]$.
\end{description}

\section{Koch snowflake}
The Koch snowflake can be constructed by starting with an equilateral triangle, then recursively altering each line segment as follows.
\par 1. Divide the line segment into three segments of equal length.
\par 2. Draw an equilateral triangle that has the middle segment from step 1 as its base and points outward.
\par 3. Remove the line segment that is the base of the triangle from step 2.
\par
\begin{figure}[hbtp]
            \centering
            \begin{tabular}{c}
              \includegraphics[width=.5\textwidth]{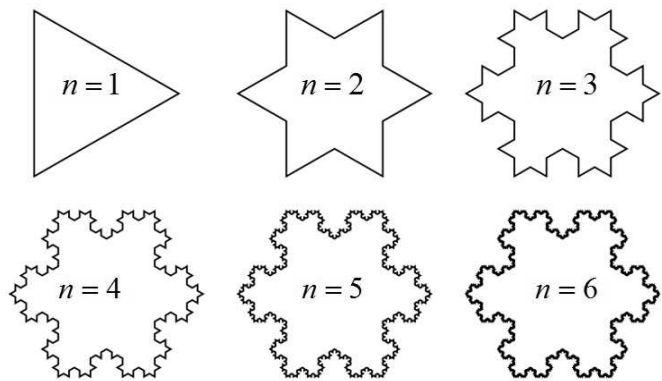}
            \end{tabular}
            \caption{The iterative process of a Koch snowflake.}
            \label{fig-sf}
 \end{figure}

\begin{figure}[hbtp]
            \centering
            \begin{tabular}{c}
              \includegraphics[width=.9\textwidth]{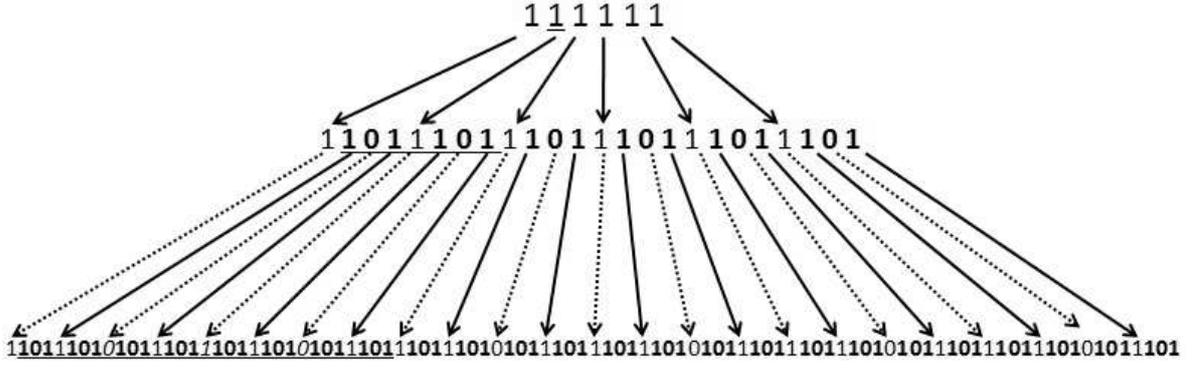}
            \end{tabular}
            \caption{The iterative process of a Koch snowflake code.}
            \label{fig-sit}
 \end{figure}

 The above process are shown as in Figure \ref{fig-sf}. After each iteration, the number of sides of the Koch snowflake increases by a factor of $3$, so the number of sides or points after $n$ iterations is given by
\begin{equation}
  3\times 4^{n-1}, \quad n=1,2,3,\cdots\cdots.
\end{equation}
Similarly as in the Koch curve, we can get the positions of the sharp points at the step $n$, where $n>1$. Hence, the sharp point will appear at
\begin{align}\label{SN-p}
\begin{split}
&4^l(4k-2)+1,\quad l=0,1,2,\cdots, n-2,\quad k=1,2,3,\cdots, 3\times 4^{n-2-l},\\
&4^{n-1}(k-1)+1, \quad k=1,2,3.
\end{split}
\end{align}
Using the same method as the above section, we get the code at the step $n$.
\begin{align}\label{SN-Cod}
\begin{split}
&n=2,\qquad 111111;\\
&n=3, \qquad 1{\bf101}1{\bf101}1{\bf101}1{\bf101}1{\bf101}1{\bf101},\\
&n=4,\\
&1{\bf101}1{\bf101}0{\bf101}1{\bf101}1{\bf101}1{\bf101}0{\bf101}1{\bf101}1{\bf101}1{\bf101}0{\bf101}1{\bf101}\\
&1{\bf101}1{\bf101}0{\bf101}1{\bf101}1{\bf101}1{\bf101}0{\bf101}1{\bf101}1{\bf101}1{\bf101}0{\bf101}1{\bf101}\\ \\
& \cdots\cdots
\end{split}
\end{align}
The iterative process can be described by Figure \ref{fig-sit}. There is an $101$ between arbitrary two elements of the previous step, which implies two sharp angle pairs and one obtuse angle pair will be inserted between arbitrary two pairs.
 \par Locally, it is the same iteration as the Koch curve, that is, each iteration will generate $4$ times of itself, which is shown in Figure \ref{fig-sit} using the underlined elements. Then we have
 \begin{prop} An affine Koch snowflake curve can be encoded as Eq. (\ref{SN-Cod}),
where the element $1$ means there are two points with affine curvatures $\kappa_1=\kappa_2=-1, \bar{\kappa}_1=-1, \bar{\kappa}_2=1$, and the element $0$ implies the two points have the affine curvatures $\kappa_1=\kappa_2=\bar{\kappa}_1=\bar{\kappa}_2=1$. Further more, from the bold number we can find the characterization. There are $6\times4^{n-2}$ elements at the step $n$, where $n\ge 2$.
\end{prop}
\begin{algorithm}[hbtp]
\label{alg-SN}
 \KwIn{initial three starting points $P_1, P_2, P_3$, and an integer $N>1$;}
 \KwOut{an affine Koch snowflake curve at the $N$th step;}
 join points $P_1$, $P_2$ and $P_3$ with two straight line segments in turn\;
 \For{$idx\Leftarrow 1$ \KwTo $6\times4^{N-2}$}{
      $i\Leftarrow 0, flag\Leftarrow 0$\;
      \If{$idx\in\{1, \frac{4^{N-1}}{2}+1,4^{N-1}+1\}$}{
        $flag\Leftarrow 1$;
      }
     \While{$i<=\log_4 (idx-1)$ {\bf and} $flag==0$}{
         \If{$\frac{idx-1}{4^i}$ is odd}{
               $flag\Leftarrow 1$;
        }
       $i\Leftarrow i+1$;
      }
     \eIf{$flag==1$}
     {use the affine curvatures $\kappa=\{-1,-1\},\bar{\kappa}=\{-1,1\} $ and Eq. (\ref{Ite-Pt}) to generate two points, start from $P_3$, join them with line segments in turn, assign into $P_1,P_2,P_3$ with the last three points\;
     }{use the affine curvatures $\kappa=\{1,1\},\bar{\kappa}=\{1,1\} $ and Eq. (\ref{Ite-Pt}) to generate two points, start from $P_3$, join them with line segments in turn, assign into $P_1,P_2,P_3$ with the last three points\;
     }
  }
\caption{Generate directly the affine Koch snowflake curve using affine curvatures}
\end{algorithm}
 \par From Eq. (\ref{SN-p}), we can locate the element $1$ at the sequence of the $n$th step,
\begin{align}
\begin{split}
&4^l(2k-1)+1,\quad l=0,1,2,\cdots, n-2,\quad k=1,2,3,\cdots, 3\times 4^{n-2-l},\\
&\frac{4^{n-1}(k-1)}{2}+1, \quad k=1,2,3.
\end{split}
\end{align}
Hence, according to the above results, by the affine curvatures, we can generate directly an affine Koch snowflake using the algorithm \ref{alg-SN} for any arbitrary step $n$. Similarly, we have
\begin{rem}
If we choose three initial points $\vec{r}_1,\vec{r}_2$ and $\vec{r}_3$ which are not collinear, the curve can be determined uniquely. If we choose $\vec{r}_3=\vec{r}_2+\left(
                       \begin{array}{cc}
                         \cos\frac{\pi}{3} & -\sin\frac{\pi}{3} \\
                         \sin\frac{\pi}{3} & \cos\frac{\pi}{3} \\
                       \end{array}
                     \right)
(\vec{r}_2-\vec{r}_1),$ it is a standard Koch snowflake curve. Otherwise, it is affinely equivalent to the standard Koch snowflake curve.
\end{rem}
\section{Hilbert curve}
\begin{figure}[hbtp]
            \centering
            \begin{tabular}{cc}
              \includegraphics[width=.4\textwidth]{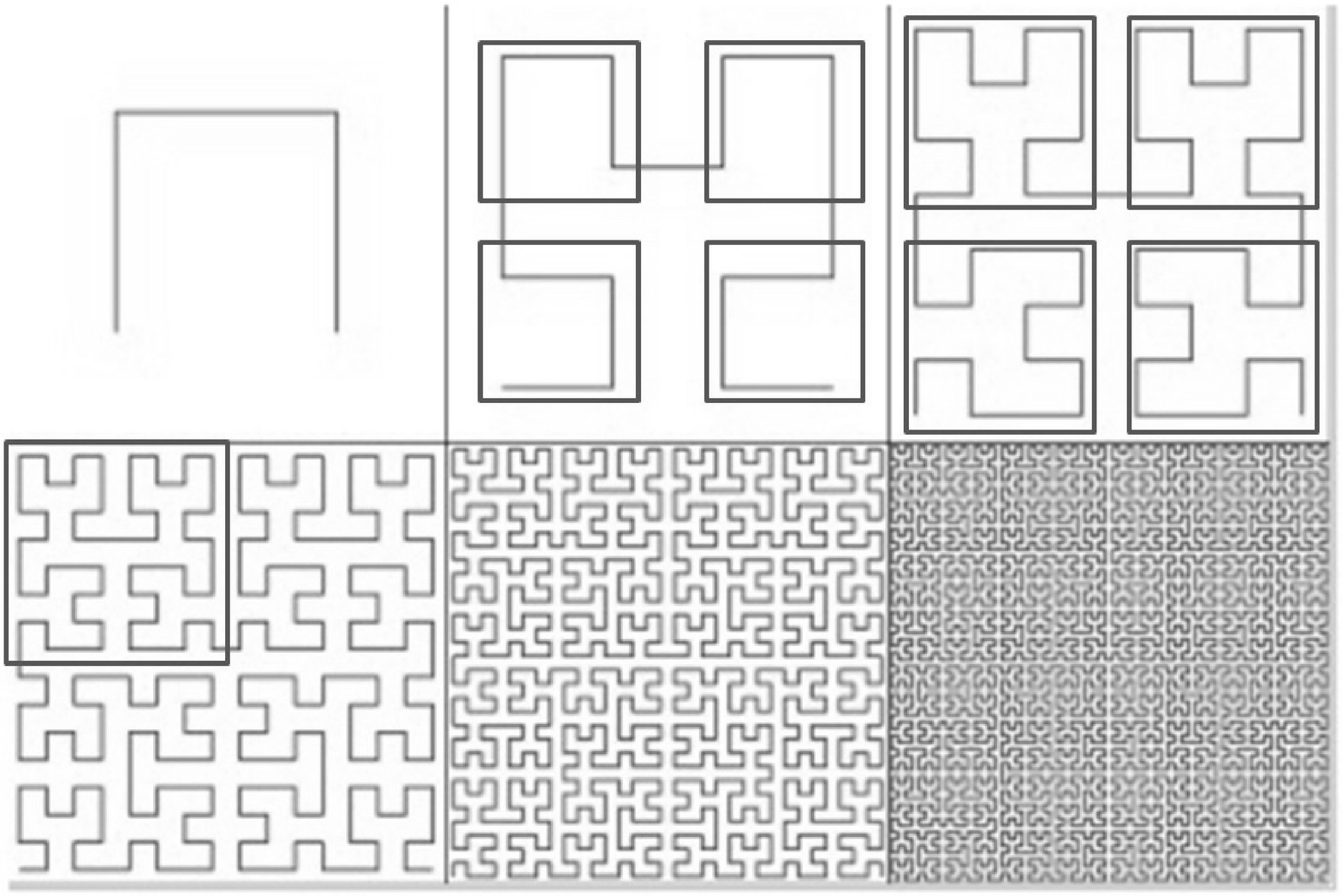}\quad&\quad\includegraphics[width=.3\textwidth]{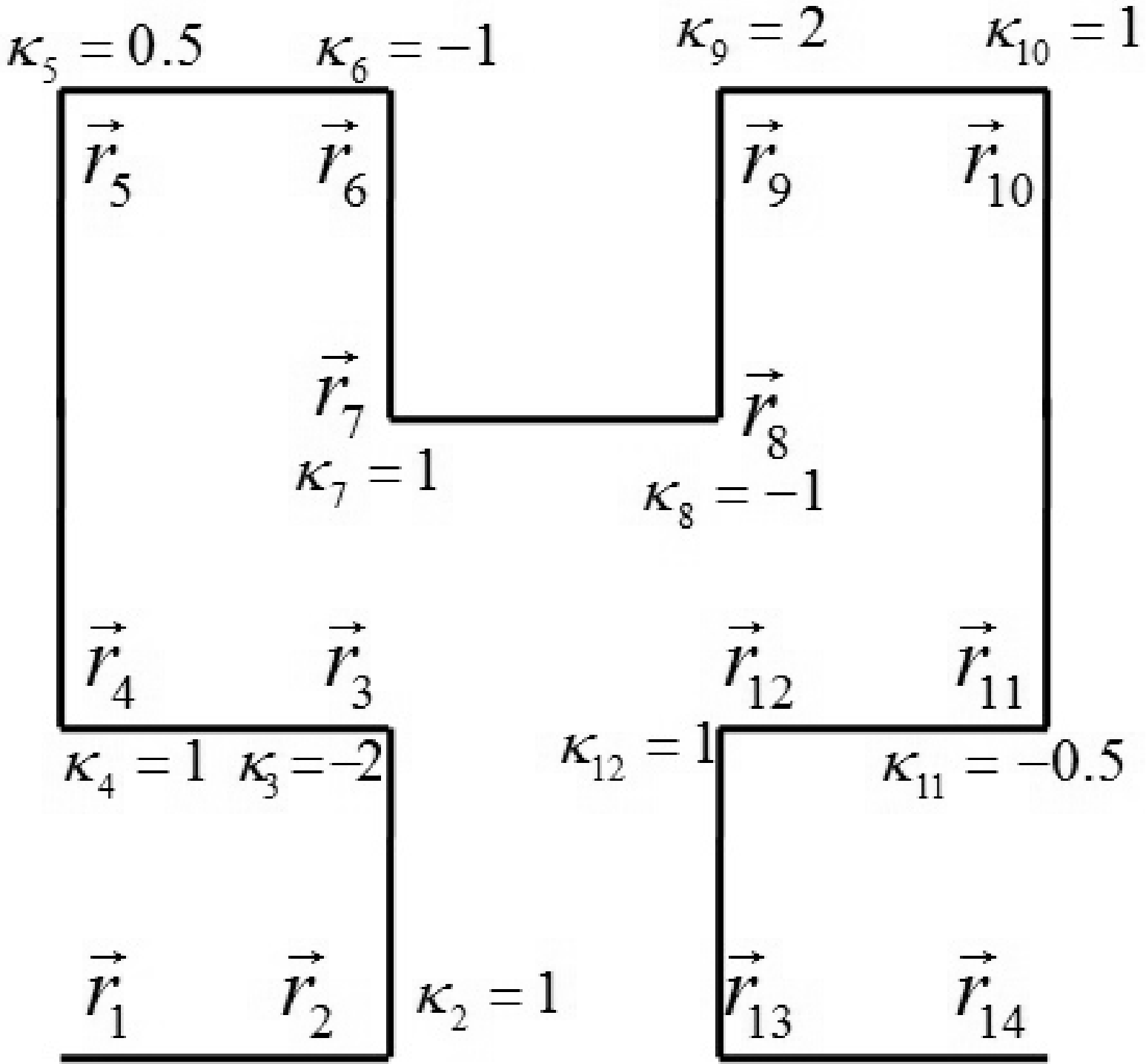}
            \end{tabular}
            \caption{{\it Left:} Six iterations of the Hilbert curve construction; {\it Right:} The inflection points on a Hilbert curve and their affine curvatures. }
            \label{fig-hil}
 \end{figure}

A Hilbert curve (also known as a Hilbert space-filling curve)  shown in the left of Figure \ref{fig-hil} is a continuous fractal space-filling curve first described by the German mathematician David Hilbert in 1891, as a variant of the space-filling Peano curves discovered by Giuseppe Peano in 1890.
 \par Now we observe the iterative process by using an affine transformation viewpoint. At every step, there are $4$ graphs in the boxes as shown in Figure \ref{fig-hil}, which are affinely equivalent to each other. Exactly, arbitrary one of the graphs in $4$ boxes at the step $n+1$ is affinely equivalent to the whole graph of the step $n$. Hence, it is important to calculate the curvatures at the joints.
 \par Let us consider the affine curvatures $\kappa$ and $\bar{\kappa}$ of Hilbert curve in details. Traditionally, there are $4^n$ points after the $n^{\mathrm{th}}$ iteration. Notice that in this paper we only use its inflection points, that is, if some points in the same segment, we just count the starting point and the end point of the segment, which is shown in the right of Figure \ref{fig-hil}.
 \par Obviously, in the Hilbert curve,  $\vec{t}_{k-1}\|\vec{t}_{k+1}, \forall k\in \mathbb{Z}^+$. By Eq. (\ref{PCur}), we have
 \begin{equation}
   \bar{\kappa}_k=0, \quad \forall k\in \mathbb{Z}^+.
 \end{equation}
 Hence, the iterative equation (\ref{Ite-Pt}) can be written as
 \begin{equation}\label{hlIte-Pt}
   \vec{r}_{k+2}=\kappa_k\vec{r}_{k-1}-\kappa_k\vec{r}_k+\vec{r}_{k+1}.
 \end{equation}
In the right of Figure \ref{fig-hil}, the first affine curvatures have been identified by a direct computation from Eq. (\ref{PCur}). In fact, the Hilbert curve is symmetry, and it also has the following properties.
\begin{lem}
A Hilbert curve $\vec{r}_1, \vec{r}_2,\cdots, \vec{r}_n$ is affinely equivalent to its inverted sequence curve $\vec{r}_n, \vec{r}_{n-1}, \cdots, \vec{r}_1$. Hence,
no matter which one of $\vec{r}_1$ and $\vec{r}_n$ is chosen as the start point, the curvatures are same at the corresponding points.
\end{lem}
\begin{figure}[hbtp]
            \centering
            \begin{tabular}{ccc}
              \includegraphics[width=.12\textwidth]{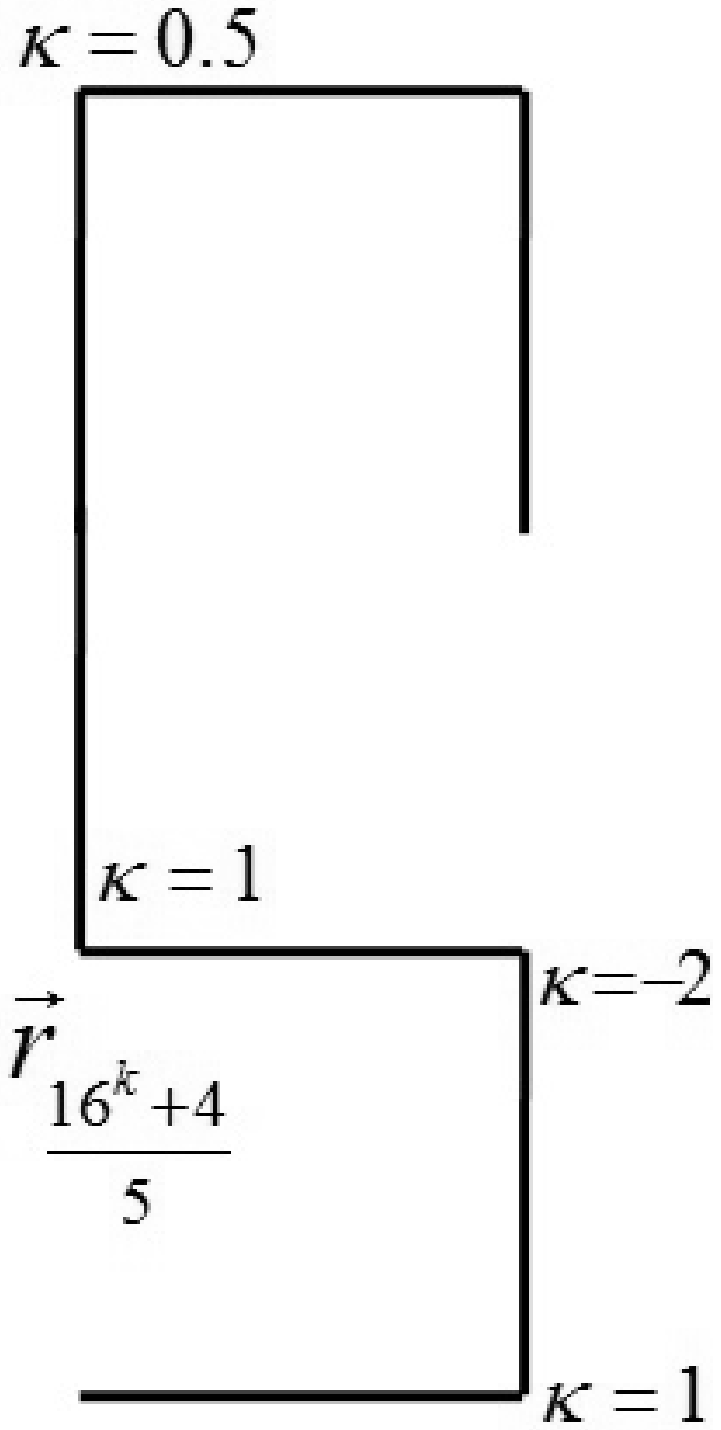}\quad\qquad&\includegraphics[width=.24\textwidth]{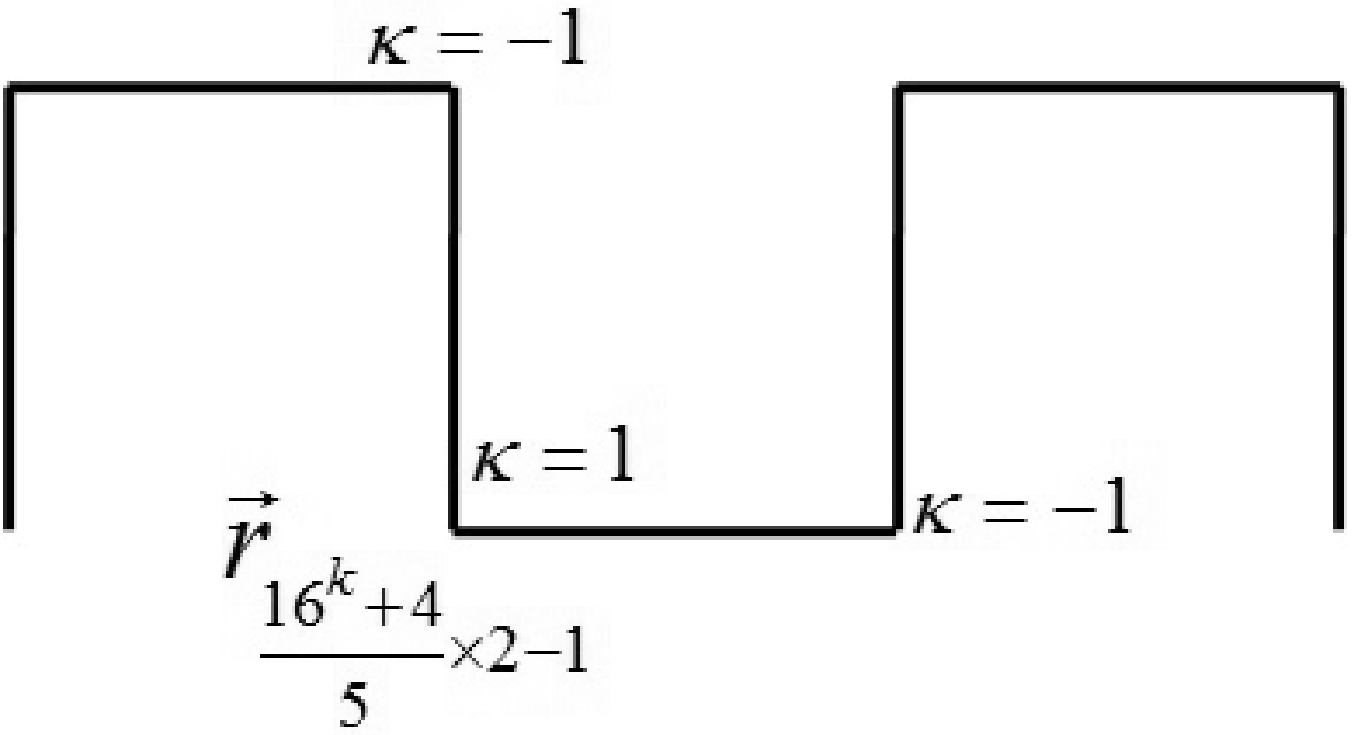}& \qquad \includegraphics[width=.16\textwidth]{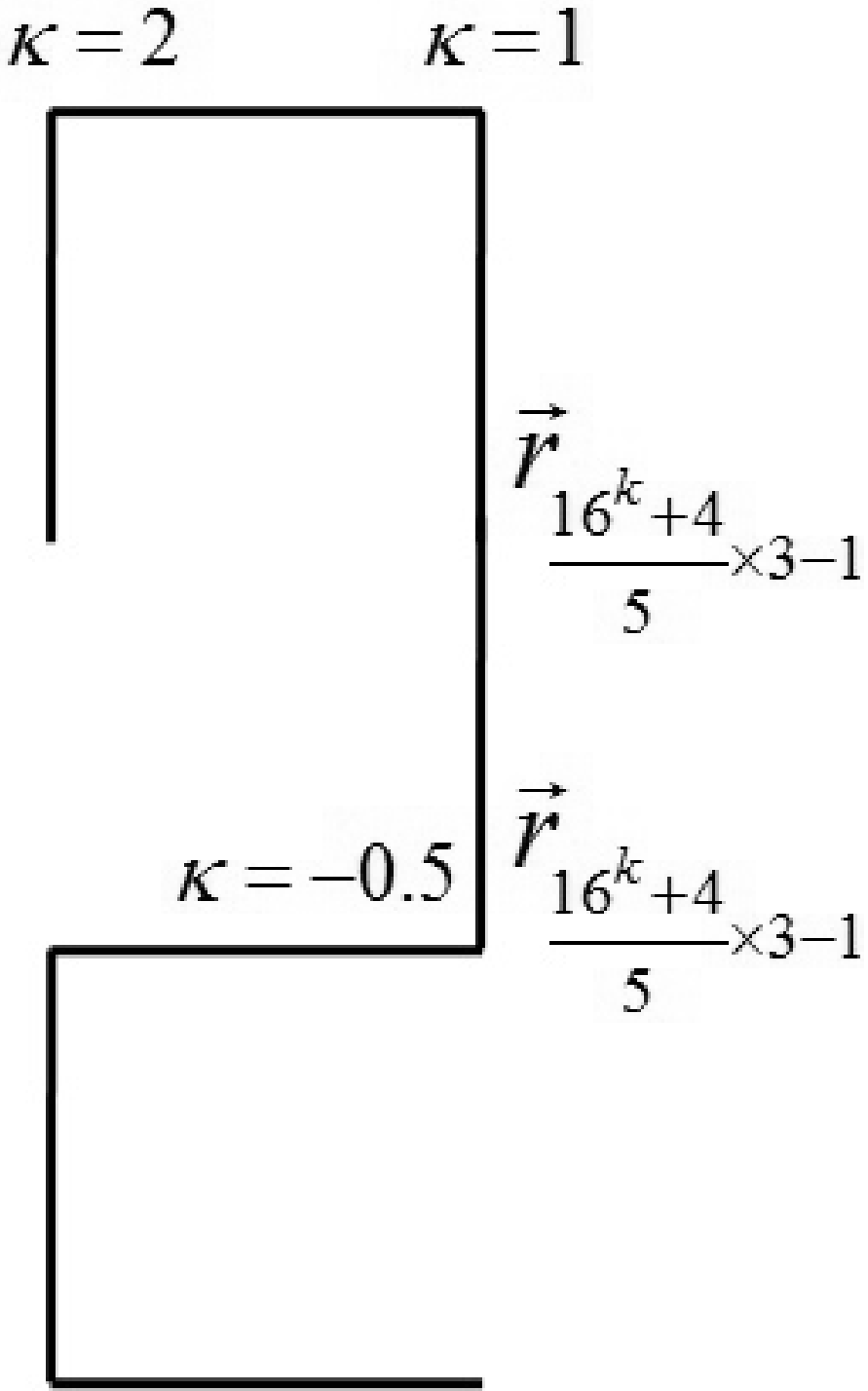}
            \end{tabular}
            \caption{The joints of iterative process from the $2k-1$st step to the $2k$th step. }
            \label{fig-o2e}
 \end{figure}
\par During the iterative process from the $2k-1$st step to the $2k$th step, at the joints,  the cases shown in Figure \ref{fig-o2e}  will appear in turn. Assume there are $N(2k-1)$ inflection points after the $2k-1$st step iteration. From the iterative process, we know the number of the inflection points after the $2k$th step iteration is
\begin{equation}\label{Num-Even}
  N(2k)=4N(2k-1)-2.
\end{equation}
 \par
\begin{figure}[hbtp]
            \centering
            \begin{tabular}{ccc}
              \includegraphics[width=.16\textwidth]{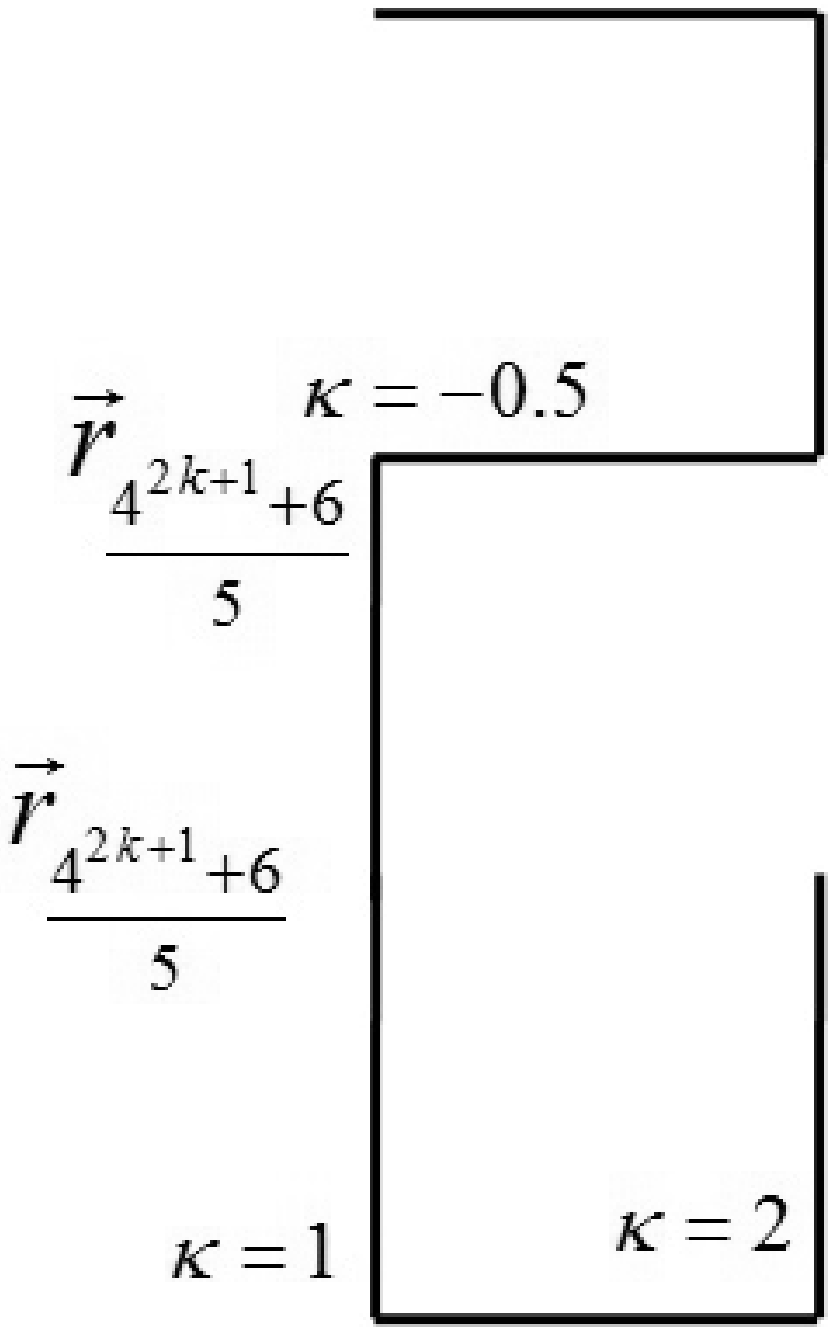}\quad\qquad&\includegraphics[width=.32\textwidth]{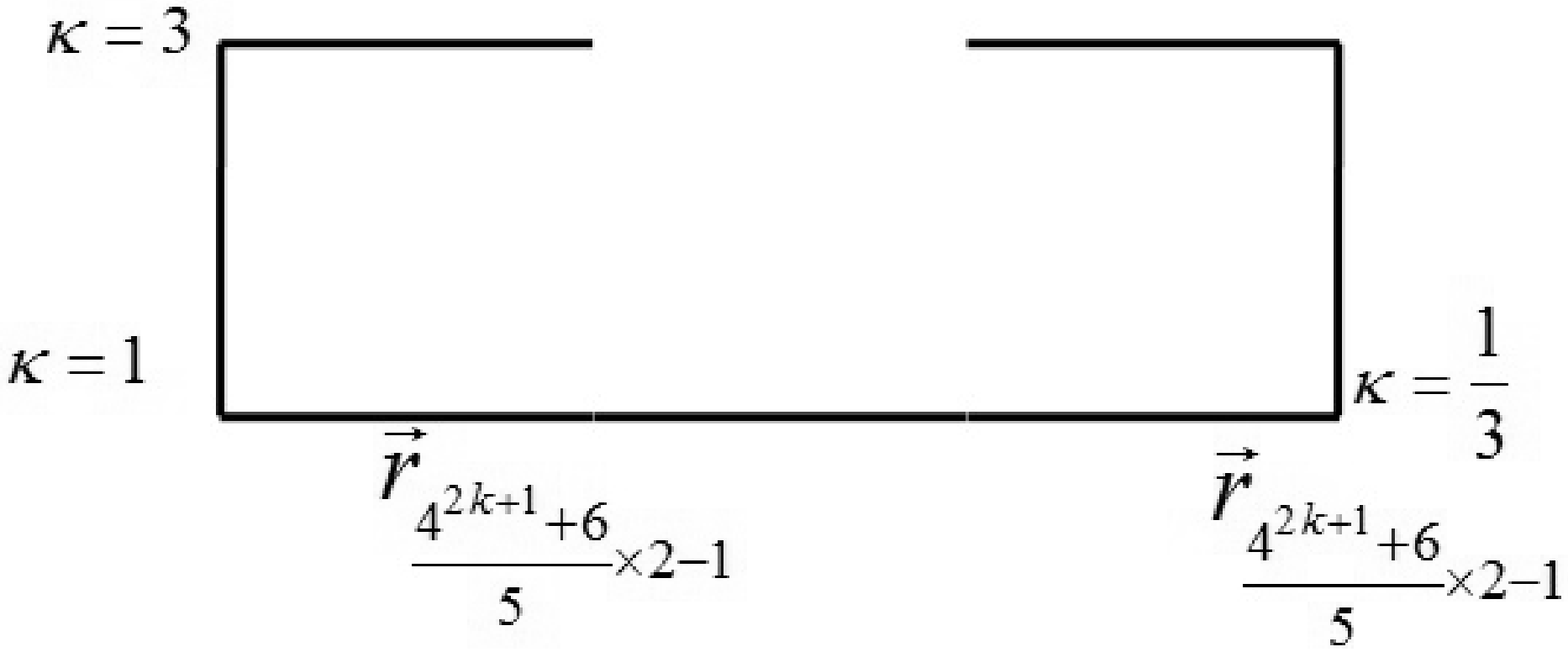}& \qquad \includegraphics[width=.20\textwidth]{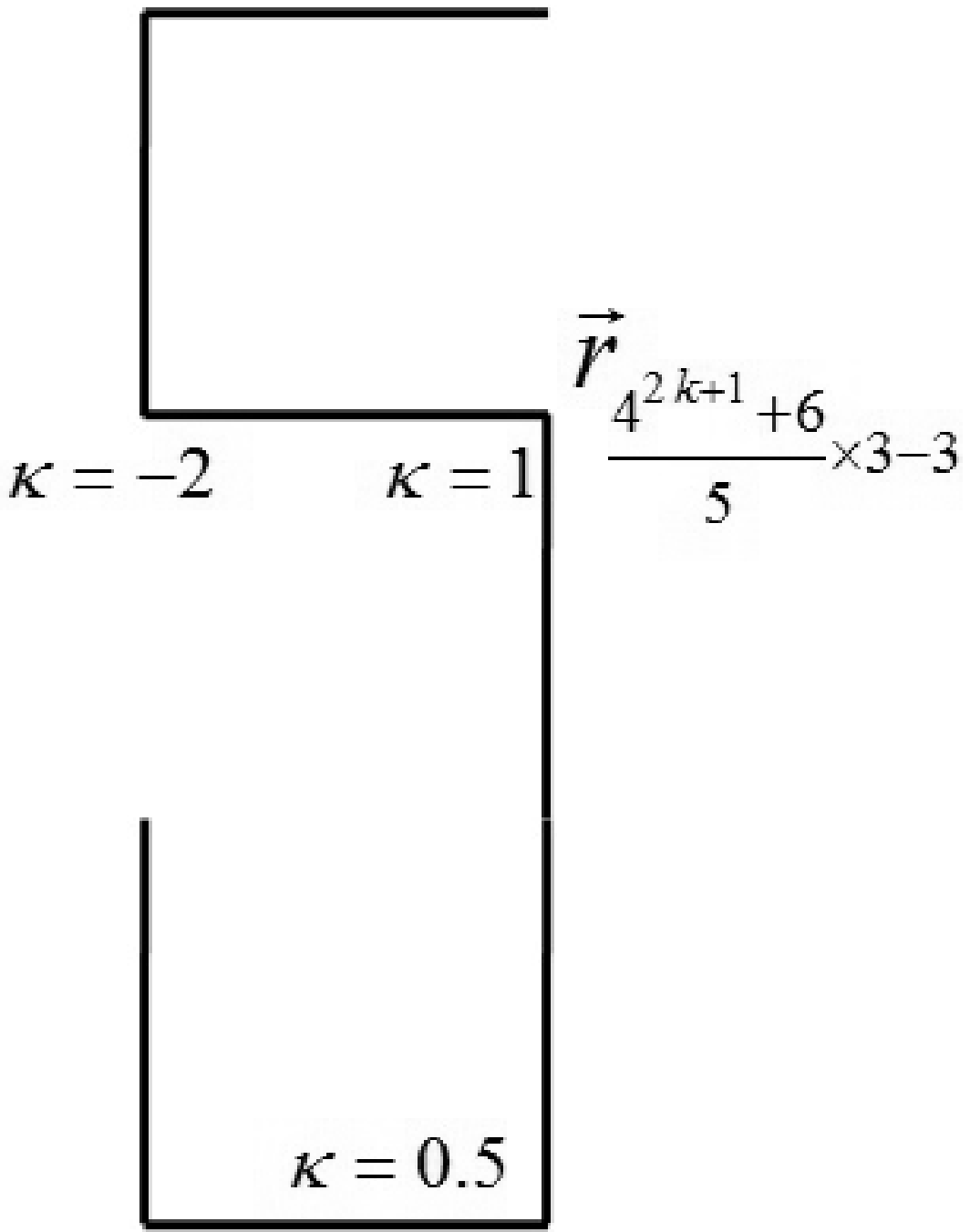}
            \end{tabular}
            \caption{The joints of iterative process from the $2k$th step to the $2k+1$st step. }
            \label{fig-e2o}
 \end{figure}

 Using the same method, let us consider the iterative process from the $2k$th step to the $2k+1$st step. At the joints,  the cases shown in Figure \ref{fig-e2o}  will appear in turn. According to these joints, it is not difficult to conclude that the number of the inflection points after the $2k+1$st step iteration is
\begin{equation}\label{Num-Odd}
  N(2k+1)=4N(2k)-4.
\end{equation}
 \par
\begin{figure}[hbtp]
            \centering
            \begin{tabular}{cc}
              \includegraphics[width=.75\textwidth]{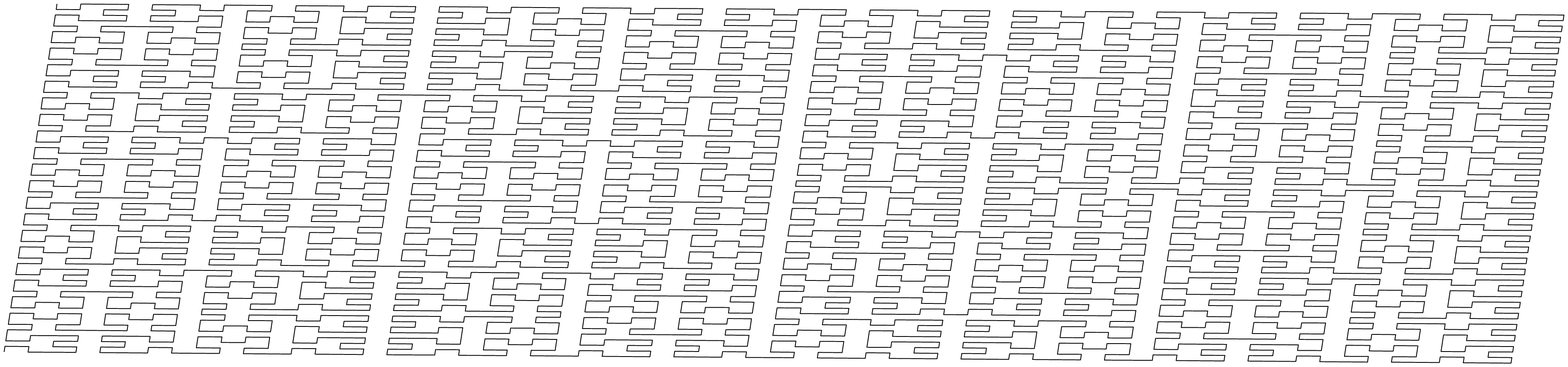}&
              \includegraphics[width=.2\textwidth]{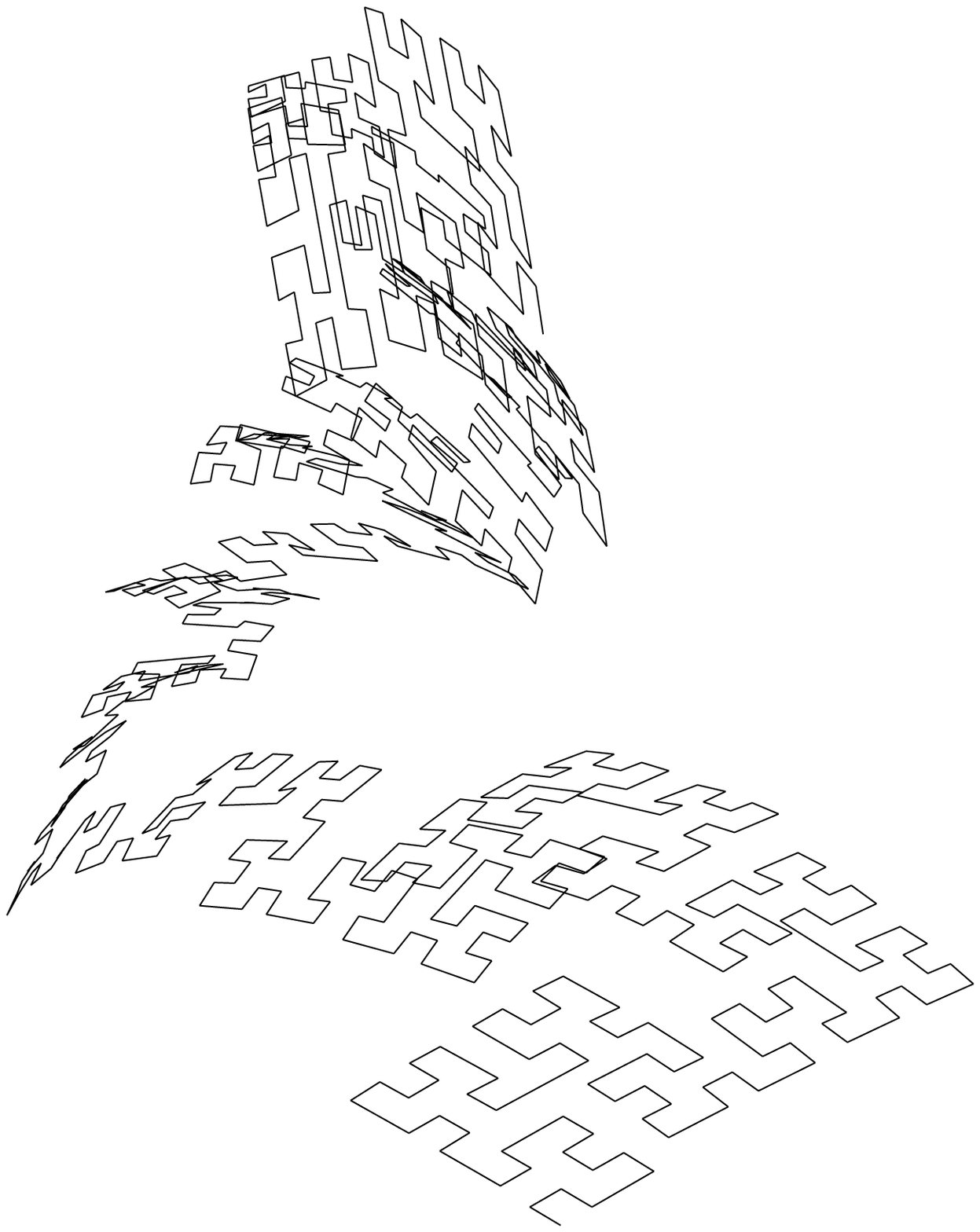}
            \end{tabular}
            \caption{Discrete affine curve generated by the affine curvatures of Hilbert curve. }
            \label{fig-HL-alg}
 \end{figure}

 \noindent Since $N(1)=4$, from Eqs. (\ref{Num-Even}) ans (\ref{Num-Odd}) we have
 \begin{lem}
 The number of the inflection points of a Hilbert curve at the $n$th step is
\begin{spacing}{1.3}
 \begin{equation}
   N(n)=\left\{
          \begin{array}{ll}
             \frac{4^{n+1}+4}{5}, & n~\hbox{is odd;}\\
           \frac{4^{n+1}+6}{5}, & n~\hbox{is even.}
          \end{array}
        \right.
 \end{equation}
\end{spacing}
 \end{lem}
Now we can calculate the change of the affine curvatures from the step $2k-1$ to the step $2k$. Exactly, after the step $2k-1$, we have obtained the curvatures
$$\kappa(2),\kappa(3),\cdots,\kappa(\frac{4^{2k}+4}{5}-2).$$
In Figure \ref{fig-o2e}, the curvatures at the inflection point have been marked. Hence we can easily write out the curvature at every inflection point in a Hilbert curve at the step $2k$.
\begin{align}
\begin{split}
&\kappa(\frac{4^{2k}+4}{5}-1)=-2,\quad \kappa(\frac{4^{2k}+4}{5})=1,\quad \kappa(\frac{4^{2k}+4}{5}+1)=0.5,\\
&\kappa(\frac{4^{2k}+4}{5}+l)=\kappa(l+1), \qquad l=2,3,\cdots,\frac{4^{2k}+4}{5}-3,\\
&\kappa(\frac{4^{2k}+4}{5}\times2-2)=-1,\quad \kappa(\frac{4^{2k}+4}{5}\times2-1)=1,\quad \kappa(\frac{4^{2k}+4}{5}\times2)=-1,\\
&\kappa(\frac{4^{2k}+4}{5}\times2+l)=\kappa(l+1), \qquad l=1,2,\cdots,\frac{4^{2k}+4}{5}-4,\\
&\kappa(\frac{4^{2k}+4}{5}\times3-3)=2,\quad \kappa(\frac{4^{2k}+4}{5}\times3-2)=1,\quad \kappa(\frac{4^{2k}+4}{5}\times3-1)=-0.5,\\
&\kappa(\frac{4^{2k}+4}{5}\times3-1+l)=\kappa(l+1), \qquad l=1,2,\cdots,\frac{4^{2k}+4}{5}-3.
\end{split}
\end{align}
Similarly, after the step $2k$, we obtain the curvatures
$$\kappa(2),\kappa(3),\cdots,\kappa(\frac{4^{2k+1}+6}{5}-2).$$
It is time to get the iterations of the affine curvatures from the step $2k$ to the step $2k+1$. From Figure \ref{fig-e2o}, we have the following equations.
Since the curvature $\kappa(\frac{4^{2k+1}+6}{5}-2)$ will be changed, we need to remember it temporarily. Taking the notation $\kappa~'=\kappa(\frac{4^{2k+1}+6}{5}-2),$ then the curvatures in turn are
\begin{align}
\begin{split}
&\kappa(\frac{4^{2k+1}+6}{5}-2)=2,\quad\kappa(\frac{4^{2k+1}+6}{5}-1)=1,\quad \kappa(\frac{4^{2k+1}+6}{5})=-0.5,\\
&\kappa(\frac{4^{2k+1}+6}{5}+l)=\kappa(l+1),\quad \qquad l=1,2,3,\cdots,\frac{4^{2k+1}+6}{5}-4,\\
&\kappa(\frac{4^{2k+1}+6}{5}\times2-3)=3,\quad \kappa(\frac{4^{2k+1}+6}{5}\times2-2)=1,\quad \kappa(\frac{4^{2k+1}+6}{5}\times2-1)=\frac{1}{3},\\
&\kappa(\frac{4^{2k+1}+6}{5}\times2-1+l)=\kappa(l+2), \quad \qquad l=1,2,\cdots,\frac{4^{2k+1}+6}{5}-5,\\
&\kappa(\frac{4^{2k+1}+6}{5}\times3-5)=\kappa~',\\
&\kappa(\frac{4^{2k+1}+6}{5}\times3-4)=-2,\quad \kappa(\frac{4^{2k+1}+6}{5}\times3-3)=1,\quad \kappa(\frac{4^{2k+1}+6}{5}\times3-2)=0.5,\\
&\kappa(\frac{4^{2k+1}+6}{5}\times3-2+l)=\kappa(l+2), \quad \qquad l=1,2,\cdots,\frac{4^{2k+1}+6}{5}-5,\\
&\kappa(\frac{4^{2k+1}+6}{5}\times4-6)=\kappa~'.
\end{split}
\end{align}
Now we conclude
\begin{prop}The affine curvatures of a Hilbert curve at every inflection point satisfy that
$$\kappa\in \{1,-1,2,-2,0.5,-0.5, 3,\frac{1}{3}\}, \quad \bar{\kappa}=0.$$
\end{prop}
Here, let us take the notation
\begin{equation}\label{A-D}
  A=-2~1~0.5,\quad B=-1~1~-1,\quad C=2~1~-0.5,\quad D=3~1~\frac{1}{3}.
\end{equation}
We can show the results of the iteration using the following style.
\begin{align}\label{ABC}
\begin{split}
& n=2,\quad 1ABC1;\\
& n=3, \quad 1ABCC1ABCDABC1AABC1.
\end{split}
\end{align}
If we use $K_3$ to represent the part between the first $1$ and the last $1$, that is
\begin{equation}
  K_3=ABCC1ABCDABC1AABC,
\end{equation}
then
\begin{equation}
n=4, \quad 1K_31AK_31B1K_3C1K_31.
\end{equation}
Hence it is easy to find the regular pattern.
\begin{align}\label{AK-ite}
\begin{split}
&1K_{2n}1= 1K_{2n-1}1AK_{2n-1}1B1K_{2n-1}C1K_{2n-1}1,\\
&1K_{2n+1}1=1K_{2n}C1K_{2n}DK_{2n}1AK_{2n}1.
\end{split}
\end{align}
The length of the sequence is
\begin{spacing}{1.3}
\begin{equation}
  y(n)=\left\{
         \begin{array}{ll}
           \frac{6\times 4^{n-1}+1}{5}, & n~\hbox{is even;} \\
           \frac{6\times 4^{n-1}-1}{5}, & n~\hbox{is odd and }n>1.
         \end{array}
       \right.
\end{equation}
\end{spacing}
From Eqs. (\ref{A-D})-(\ref{AK-ite}), if we take the notation
\begin{gather}\label{P-nota}
\begin{split}
  P=ABC=-2~1~0.5~-1~1~-1~2~1~-0.5,\quad S=C1=2~1~-0.5~1,\\
  T=D= 3~1~\frac{1}{3},\quad U=1A=1~-2~1~0.5,\quad V=1B1=1~-1~1~-1~1,
\end{split}
\end{gather}
and temporarily ignore the two $1$'s at both sides of every sequence in Eq. (\ref{AK-ite}), the iterative process can be represented by
\begin{align}\label{P-ite}
\begin{split}
&n=2,\qquad  K_2=P;\\
&n=3,\qquad  K_3=PSPTPUP;\\
&n=4,\qquad  K_4=K_3UK_3VK_3SK_3;\\
&n=5,\qquad  K_5=K_4SK_4TK_4UK_4;\\
&n=6,\qquad  K_6=K_5UK_5VK_5SK_5;\\
&\cdots\quad\cdots\quad\cdots
\end{split}
\end{align}
By observing the sequences, we can get its regularities, which yield the following proposition.
\begin{prop} For an affinely Hilbert curve, it has the following properties:
\begin{enumerate}
  \item[I.] the affine curvatures of the affine Hilbert curve can be obtain by Eqs. (\ref{P-nota}) and (\ref{P-ite}), and the two $1$'s should be added at both side of the sequence automatically;
  \item[II.]in the sequence of the affine Hilbert curve with the step $n$,
    \begin{enumerate}
      \item $P$ will appear at the index: $2k-1, \quad k=1,\cdots, 4^{n-2}$;
      \item $T$ will appear at the index: $4^{2l-1}(2k-1), l=1,2,\cdots,\lfloor\frac{n-1}{2}\rfloor,$ $k=1,2,\cdots, 4^{n-2l-1}$, where $\lfloor x\rfloor$ is the largest integer less than or equal to $x$ ;
      \item $V$ will appear at the index: $4^{2l}(2k-1), l=1,2,\cdots,\lfloor\frac{n-2}{2}\rfloor, k=1,2,\cdots, 4^{n-2l-2}$;
      \item $S$ will appear at the index: $4^{2l-2}(8k-6),$ $l=1,2,\cdots,\lfloor\frac{n-1}{2}\rfloor,$ $k=1,2,\cdots, 4^{n-2l-1}$,
            or $4^{2l-1}(8k-2),$ $l=1,2,\cdots,\lfloor\frac{n-2}{2}\rfloor, k=1,2,\cdots, 4^{n-2l-2}$;
      \item $U$ will appear at the index: $4^{2l-2}(8k-2),$ $l=1,2,\cdots,\lfloor\frac{n-1}{2}\rfloor,$ $k=1,2,\cdots, 4^{n-2l-1}$,
            or $4^{2l-1}(8k-6),$ $l=1,2,\cdots,\lfloor\frac{n-2}{2}\rfloor, k=1,2,\cdots, 4^{n-2l-2}$;
    \end{enumerate}
  \item[III.]in the sequence of the affine Hilbert curve with the step $n$, there are $2\times4^{n-2}-1$ letters; if $n=2m$, there are $4^{2m-2}$ letter $P$'s, $\frac{4^{2m-2}-1}{15}$ letter $V$'s, $\frac{4^{2m-1}-4}{15}$ letter $T$'s, $\frac{4^{2m-2}-1}{3}$ letter $S$'s and $\frac{4^{2m-2}-1}{3}$ letter $U$'s; if $n=2m+1$, there are $4^{2m-1}$ letter $P$'s, $\frac{4^{2m-1}-4}{15}$ letter $V$'s, $\frac{4^{2m}-1}{15}$ letter $T$'s, $\frac{4^{2m-1}-1}{3}$ letter $S$'s and $\frac{4^{2m-1}-1}{3}$ letter $U$'s.
\end{enumerate}
\end{prop}
Now, the above proposition can lead to the algorithm \ref{alg-HL}. By using the algorithm \ref{alg-HL}, we could generate an affine Hilbert curve.
\begin{rem}
If we choose three initial points $\vec{r}_1,\vec{r}_2$ and $\vec{r}_3$ which are not collinear, the curve can be determined uniquely. If we choose $\vec{r}_3=\vec{r}_2+\left(
                       \begin{array}{cc}
                         0 & -1 \\
                         1 & 0 \\
                       \end{array}
                     \right)
(\vec{r}_2-\vec{r}_1),$ it is a standard Hilbert curve. Otherwise, it is affinely equivalent to the standard Hilbert curve.
\end{rem}
\begin{algorithm}[hbtp]
\label{alg-HL}
 \KwIn{initial three starting points $P_1, P_2, P_3$, and an integer $N>1$;}
 \KwOut{an affinely Hilbert curve at the $N$th step;}
 join points $P_1$ $P_2$ and $P_3$ with two straight line segments in turn\;
 use $\kappa=1$ and Eq. (\ref{hlIte-Pt}) to obtain the fourth point $P$,  join it to point $P_3$ with a straight line segment\;
            $P_1\Leftarrow P_2, P_2\Leftarrow P_3, P_3\Leftarrow P$\;
 \For{$idx\Leftarrow 1$ \KwTo $2\times4^{N-2}-1$}{
      $i\Leftarrow 0, flag\Leftarrow 0$\;
      \If{$idx$ is odd}{
        $flag\Leftarrow 'P'$;
      }
     \While{$i<=\frac{\log_4 (idx-1)}{2}+1$ {\bf and} $flag==0$}{
            \uIf{$\mod(idx/4^{2i}, 2)==1$}{
            $flag\Leftarrow 'V'$;
            }
            \uElseIf{$\mod(idx/4^{2i-1}, 2)==1$}{
            $flag\Leftarrow 'T'$;
            }
            \uElseIf{$\mod(idx/4^{2i-2}, 8)==2$ {\bf or} $\mod(idx/4^{2i-1}, 8)==6$}{
            $flag\Leftarrow 'S'$;
            }
            \Else{
            $flag\Leftarrow 'U'$;
            }
       $i\Leftarrow i+1$;
      }
      use Eqs. (\ref{hlIte-Pt}), (\ref{P-nota}) and $flag$ to generate the corresponding points in turn, start from $P_3$, join them with line segments in turn, assign into $P_1,P_2,P_3$ with the last three points\;
  }
     use $\kappa=1$ and Eq. (\ref{hlIte-Pt}) to obtain the last point,  join it to point $P_3$ with a line segment.
\caption{Generate an affinely Hilbert curve using affine curvatures}
\end{algorithm}
The following example shows how to use the affine curvatures of the Hilbert curve to generate fractal curves.\\
\noindent {\bf Example 2.} In this example, we give some curves which are obtained by the affine curvatures of the affine Hilbert curve:
\begin{description}
  \item Firstly, using the initial points $[0~0;-1~-2; -10~1]$ and $n=6$, we obtain an affine Hilbert curve as shown in the first one of Figure \ref{fig-HL-alg}.
  \item Then, using the the initial points $[0~0~1;0~1~1;1~1~1]$ and $n=7$, with the same affine curvatures as the affine Hilbert curve, and $\bar{\kappa}_i=\tau_i=0.005\kappa_i$, we get a discrete space curve as shown in the second graph of Figure \ref{fig-HL-alg}.
\end{description}
\section{Conclusions}
By the affine curvatures, the Koch curves can be encoded as one of the following forms $$1, 1{\bf 0}1{\bf1}1{\bf0}1, 1011101{\bf0}1011101{\bf1}1011101{\bf0}1011101, \cdots\cdots,$$
where the element $1$ means there are two points with affine curvatures $\kappa_1=\kappa_2=-1, \bar{\kappa}_1=-1, \bar{\kappa}_2=1$, and the element $0$ implies the two points have the affine curvatures $\kappa_1=\kappa_2=\bar{\kappa}_1=\bar{\kappa}_2=1$. Furthermore, from the bold number we can find the characterization of the iteration. There are $2\times4^{n-2}-1$ elements at the step $n$, where $n\ge 2$.
\par The affine curvatures of the affine Hilbert curve can be obtain by Eqs. (\ref{P-nota}) and (\ref{P-ite}), and the two $1$'s should be added at both side of the sequence automatically.
\par By these codes, it is easy, quick and direct to generate the corresponding fractal curves. Of course, the regularities are more obvious. In this paper, as a start, we consider the Koch curve, the Koch snowflake curve and the Hilbert curve. However, there are still many discrete fractal curves needed to be studied for future works. After we obtain regularities of their curvatures, the fractal curves can be classified by the affine equivalence property. Since the fractal curves have been quantified by the affine curvatures, using a fast fourier transform(FFT) algorithm, the periodicity can be find.
\renewcommand{\refname}{\bf\fontsize{12}{12}\selectfont References}
\bibliographystyle{amsplain}

\end{document}